\documentclass[a4paper,12pt,reqno]{amsart}
\usepackage{amssymb,amsmath,amsfonts,latexsym} %%
\usepackage{hyperref}
\usepackage{ifthen}
\usepackage{graphicx}
\usepackage{float}
\usepackage{mathrsfs}
\usepackage{comment}
\usepackage{enumerate}
\usepackage[foot]{amsaddr}
 \usepackage[left=25mm,right=25mm,top=25mm,bottom=25mm,paper=a4paper]{geometry}

\theoremstyle{plain}

\newtheorem*{thm*}{Theorem}
\numberwithin{equation}{section}

\theoremstyle{definition}
\newcounter {own}
\def\theown {\thesection  .\arabic{own}}

{\qed\bigskip}

\newcounter{alphabet}

\newcommand{\R}{\mathbb{R}}

\newcommand{\T}{\mathcal{T}}

\newtheorem{theorem}{Theorem}[section]

\newtheorem{definition}{Definition}[section]
\newtheorem{remark}{Remark}[section]
\newtheorem{lemma}{Lemma}[section]
\newtheorem{prop}{Proposition}[section]
\newcounter{minutes}
\divide\time by 60
\newcounter{hours}
\multiply\time by 60 \addtocounter{minutes}{-\time}
\begin{document}
 
 %%% Pls. specify the kind of the article:
 %%% This is: %% ORIGINAL PAPER %%%% 
 %%% or: REVIEW PAPER / SHORT PAPER, etc. %%%%%%%%

%%%% Title of article for FCAA %%%%%%%%%%%%%%%%%%%%%%%%
\title{ on the boundedness of generalized fractional integral operators in Morrey spaces and Campanato spaces associated with the Dunkl operator on the real line}

%=========================================================================
\thanks{%$^\dagger$
File:~\jobname .tex,
          printed: \number\year-\number\month-\number\day,
          \thehours.\ifnum\theminutes<10{0}\fi\theminutes}
%=========================================================================

% \titlerunning{ON THE BOUNDEDNESS OF A GENERALIZED DUNKL-TYPE FRACTIONAL
% INTEGRAL}
%% if too long for running head - use the text from 1st line !%

%%%%% authors:
\author{Sumit Parashar}
 \address{Sumit Parashar, School of Mathematical and Statistical Sciences, Indian Institute of Technology Mandi, Kamand,
 Mandi, HP 175005, India \newline
Email: d22035@students.iitmandi.ac.in }
% \email{d22035@students.iitmandi.ac.in} 

\author{ Saswata Adhikari $^\dagger$}
\address{Saswata Adhikari, School of Mathematical and Statistical Sciences, Indian Institute of Technology Mandi, Kamand,
 Mandi, HP 175005, India. \newline
 Email: saswata@iitmandi.ac.in  \newline
 $^\dagger$ {\tt Corresponding author}}

% %%% Dates:
% \date{Received: .... / Revised: .... / Accepted: ......}
% These dates will be entered by the Editor (EiC, V. Kiryakova) or authors/ system %

%%%%% For Production Dept.: variants of COPYRIGHT notice to appear: %%%%%%
% if Open Access option chosen, it is as:
   %% "\copyright The Author(s) 2022" % (or next year..) %
% if No Open Access, the TCA form signed to FCAA journal is used, and it appears:
   %% "\copyright Diogenes Co. Ltd. 2022" %% or 2021, or next year %%%%%
%%%%%%%%%%%%%%%%%%%%%%%%%%%%%%%%%%%%%%%%%%%%%%%%%%%%%%%%%%%%%%%%%%%%%%%%%%

%%%%%%%%%%%%%%%%%%%%%%%%%%%%%%%%%%%%%%%%%%%%%%%%%%%%%%
\begin{abstract}
It is known that the Dunkl-type fractional integral operator $I_\beta$ $(0 < \beta < 2\alpha + 2 =d_\alpha)$ is bounded from $L^p(\R,d\mu_\alpha)$ to $L^q (\R, d\mu_\alpha)$ when $1 < p < \frac{d_\alpha}{\beta}$ and $\frac{1}{p} - \frac{1}{q} = \frac{\beta}{d_\alpha}$. In \cite{spsa}, the authors introduced the generalized Dunkl-type fractional integral operator  $T_\rho^\alpha$ and it's modified version $\tilde{T}_\rho^\alpha$. Indeed they proved the boundedness of $T_\rho^\alpha$ from the generalized Dunkl-type Morrey space $L^{p,\phi}(\mathbb{R},d\mu_{\alpha})$ to the generalized Dunkl-type Morrey space $L^{q,\psi}(\mathbb{R},d\mu_{\alpha})$ for $1<p<q<\infty$ and bounedness of $\tilde{T}_\rho^\alpha$ from the generalized Dunkl-type Campanato space $\mathcal{L}^{1,\phi}(\R,d\mu_\alpha)$ to the generalized Dunkl-type Campanato space $\mathcal{L}^{1,\psi}(\R, d\mu_\alpha)$. In this paper, we study the boundedness of these operators mainly on generalized Dunkl-type Campanato spaces. However, first we prove that $ T_\rho^\alpha $ is bounded from $L^{p,\phi}(\mathbb{R},d\mu_{\alpha})$ to $L^{q,\psi}(\mathbb{R},d\mu_{\alpha})$ when $p=q=1$ and when $1<p=q<\infty$ by using Fubini’s theorem and the generalized translation on the Bessel-Kingman hypergroups. We also show that $ \tilde{T}_\rho^\alpha $ is bounded from the generalized Dunkl-type Morrey space $L^{p,\phi}(\mathbb{R},d\mu_{\alpha})$ to the generalized Dunkl-type Campanato space $\mathcal{L}^{p,\psi}(\mathbb{R},d\mu_{\alpha})$ when $1\leq p<\infty$. Further, for the case $1 < p < \infty $, we prove the boundedness of $ \tilde{T}_\rho^\alpha $ not only from the generalized Dunkl-type Morrey space to the generalized Dunkl-type Campanato space, but also between two generalized Dunkl-type Campanato spaces.
\end{abstract}

\maketitle
\pagestyle{myheadings}
\markboth{Sumit Parashar and Saswata Adhikari}{ON THE BOUNDEDNESS OF A GENERALIZED DUNKL-TYPE FRACTIONAL INTEGRAL}
 \\

\subjclass \textbf{Mathematics Subject Classification[2010]} {Primary  42B10, 47G10; Secondary 42B20, 42B25}

%%% These are examples only. Pls. use MSC 2020 for suitable topic numbers %%%%%%

 %%%%%%%%%%%%%

%%%%%%%% begin papers' body %%%%%%%%%%%%%%%%%%%%%%%%%%%%%
\section{\bf{Introduction}} \label{sec:intro}

Let $\R^+ :=(0, \infty)$. Associated to a function $\rho : \mathbb{R}^+ \to \mathbb{R}^+$, the generalized fractional integral operator $T_\rho$ is defined by 
\begin{align*}
    T_\rho f(x) := \int_{\mathbb{R}^n} f(y) \frac{\rho(|x -y|)}{|x -y|^n} dy,
\end{align*}
and its modified version $\tilde{T}_\rho$ is defined by 
\begin{align*}
\tilde{T}_\rho f(x) := \int_{\mathbb{R}^n} f(y) \left( \frac{\rho(|x -y|)}{|x-y|^n} - \frac{\rho(|y|)}{|y|^n}(1 - \chi_{B_0}(y) ) \right) dy,
\end{align*}
for any suitable function $f$ on $\R^n$.
In the above, $B_0$ and $\chi_{B_0}$ denote the unit ball around the origin and the characteristic function of $B_0$, respectively. If $\rho(t) = t^\beta , 0 < \beta < n$, then $T_\rho$ reduces to the fractional integral operator.

Next, for $1 \leq p < \infty$ and a function $\phi : \mathbb{R}^+ \to \mathbb{R}^+$, the generalized Morrey space on $\mathbb{R}^n$ is defined by 
\begin{align*}
    L^{p, \phi}(\mathbb{R}^n) := \{ f \in L^p_{loc}(\mathbb{R}^n) : ||f||_{L^{p, \phi}} < \infty \},
\end{align*}
where $||f||_{L^{p, \phi}} := \displaystyle \sup_{\substack{{r>0}\\{a \in \mathbb{R}^n}}} \frac{1}{\phi(r)} \left ( \frac{1}{r^n} \int_{B(a,r)} |f(x)|^p dx \right)^{\frac{1}{p}} $. When $\phi(r) = r^{- \frac{n}{q}}$, where $ 1 \leq p \leq q < \infty$, then the generalized Morrey space $L^{p, \phi}(\mathbb{R}^n)$ reduces to the  Morrey space $L^{p, q}(\mathbb{R}^n)$ and moreover for $p=q$, it is the Lebesgue space $L^p(\mathbb{R}^n)$.

% For $L^{p,\phi}(\R^n)$, the function $\phi(r)$ is usually satisfied to be doubling condition and $r^n \phi^p(r)$ to be almost decreasing.\\
Now, for $1 \leq p < \infty$ and a function $\phi : \mathbb{R}^+ \to \mathbb{R}^+$, the generalized Campanato space denoted by $\mathcal{L}^{p,\phi}(\mathbb{R}^n)$, is defined as
\begin{align*}
    \mathcal{L}^{p,\phi}(\mathbb{R}^n) := \{ f \in L^1_{loc}(\mathbb{R}^n) : ||f||_{\mathcal{L}^{p,\phi}} < \infty \},
\end{align*}
where $$ ||f||_{\mathcal{L}^{p,\phi}} := \displaystyle \sup_{\substack{{r>0}\\{a \in \mathbb{R}^n}}} \frac{1}{\phi(r)} \left(\frac{1}{r^n} \int_{B(a, r)} | f(x) - f_{B(a,r)}|^p dx \right)^{1/p} $$  with $f_{B(a,r)} = \frac{1}{r^n} \int_{B(a, r)} f(x) dx$, where $B(a,r)$ denotes the ball center with $a$ and radius $r$. For $p =1$, $\mathcal{L}^{1, \phi}(\R^n)$ is the $BMO_\phi(\mathbb{R}^n)$ space. In order to define $L^{p, \phi}(\R^n)$ and  $\mathcal{L}^{p,\phi}(\R^n)$, we assume that $\phi$ satisfies doubling condition and $r^n \phi^p(r)$ is almost decreasing.

Hardy-Littlewood \cite {hardy1928some,hardy1932some} and Sobolev \cite {sobolev1938theorem} proved the boundedness of the fractional integral operator on the Lebesgue space.
Further, the boundedness of fractional integral operator was extended to the Morrey space in \cite{peetre1969theory, adams1975note, chiarenza1987morrey}. The fractional integral operator $I_\beta$ is an important tool in harmonic analysis, with applications in partial differential equations. Its kernel, $|x|^{\beta -n}$, exhibits singularities at both zero and infinity. Although one might expect the generalized kernel $ \frac{\rho(|y|)}{|y|^n} $ to allow more singularity, especially involving logarithmic factors, this is not feasible within standard Lebesgue spaces. To address this, researchers have developed the theory of the generalized fractional integral operator $ T_\rho $ in broader function spaces that extend beyond Lebesgue spaces, focusing on establishing boundedness conditions for such operators. Nakai \cite{nakai2001generalized, nakai2002} showed that the generalized fractional integral operator $T_\rho$ is bounded from $L^{1, \phi}(\R^n)$ to $L^{1, \psi}(\R^n)$ and it's modified version $\tilde{T}_\rho$ is bounded from $\mathcal{L}^{1, \phi}(\R^n)$ to $\mathcal{L}^{1, \psi}(\R^n)$, under some appropriate conditions on $\rho$, $\phi$ and $\psi$. Eridani \cite{eridani2002boundedness} showed that, for $1 < p < \infty$, $T_\rho$ is bounded from $L^{p, \phi}(\R^n)$ to $L^{p, \psi}(\R^n)$ and $\tilde{T}_\rho$ is bounded from $L^{p, \phi}(\R^n)$ to $\mathcal{L}^{p, \psi}(\R^n)$, under similar conditions on $\rho$, $\phi$ and $\psi$. Further Eridani et al. \cite{eridani2004generalized} showed the boundedness of $T_{\rho}$ from $L^{p, \phi}(\R^n)$ to $L^{p, \psi}(\R^n)$ and the boundedness of $\tilde{T}_\rho$ from $\mathcal{L}^{p, \phi}(\R^n)$ to $\mathcal{L}^{p, \psi}(\R^n)$ for $1 < p < \infty$.

In $1989$, the Dunkl operators were introduced by C.F. Dunkl \cite{dunkl1989differential}, which are differential difference operators on the real line. For a real parameter $\alpha \geq -\frac{1}{2}$, the Dunkl operators are denoted by $\Lambda_\alpha$ and these are associated with the reflection group $Z_2$ on $\mathbb{R}$.  We also refer \cite{de1993dunkl, sy, rsl2} for more details on Dunkl theory. Using the Dunkl kernel, Dunkl defined the Dunkl transform $\mathcal{F}_\alpha$ in \cite{dunkl1992hankel}. R\"{o}sler in \cite{rosler1994bessel} showed that the Dunkl kernel verifies a product formula. This allows one to define the Dunkl translation $\tau_x^\alpha$, $x \in \mathbb{R}$ and as a result one has the Dunkl convolution.

For clarity, we shall now mention certain definitions and terms that are formally stated and elaborated upon in the later sections of this paper, for instant we refer \eqref{eq: 3.1a} for Dunkl-type Morrey spaces $L^{p, q}(\mathbb{R}, d\mu_\alpha)$, \eqref{eq:3.7} for generalized Dunkl-type Morrey spaces $L^{p, \phi}(\mathbb{R}, d\mu_\alpha)$, \eqref{eq: 7.1} for Dunkl-type $BMO_\phi$ spaces and \eqref{eq:3.8} for generalized Dunkl-type Campanato spaces, as well as \eqref{eq: 3.4a} for the Dunkl-type Hardy-Littlewood maximal operator $M^\alpha$, \eqref{eq:3.1} for generalized Dunkl-type fractional integral operators $T_\rho^\alpha$ and \eqref{eq:3.4} for the modified version of the generalized Dunkl-type fractional integral operators $\tilde{T}_\rho^\alpha$. 

The maximal operator and the fractional maximal operator were studied associated with Dunkl operator on $\mathbb{R}$ in \cite{abdelkefi2007dunkl,thangavelu2007riesz} and \cite{guliyev2010p}, respectively. E. Guliyev et al. \cite{guliyev2010necessary} studied  the boundedness of the Dunkl-type fractional maximal operator in the Dunkl-type Morrey space. In \cite{guliyev2009fractional}, Mammadov obtained necessary and sufficient conditions on the parameters for the
boundedness of the Dunkl-type maximal operator and the Dunkl-type
fractional integral operator from the spaces $L^{p}(\mathbb{R}, d\mu_\alpha)$ to the spaces $L^{q}(\mathbb{R},d\mu_\alpha)$, $1 < p < q < \infty$, and 
proved that the Dunkl-type modified fractional integral operator is bounded from the space $L^p(\mathbb{R}, d\mu_\alpha)$ to the Dunkl-type $BMO$ space i.e., $BMO(\mathbb{R}, d\mu_\alpha)$. In \cite{spsa}, we proved that, for $1 < p < q < \infty$,  $T_\rho^\alpha$ is bounded from $L^{p, \phi}(\R, d\mu_\alpha)$ to $L^{q, \psi}(\R, d\mu_\alpha)$, while $\tilde{T}_\rho^\alpha$ is bounded from $BMO_{\phi}(\R, d\mu_\alpha)$ to $BMO_{\psi}(\R, d\mu_\alpha)$, under some appropriate conditions on $\rho$, $\phi$ and $\psi$. However, the approach used in that proof does not follow to the case when $p = q = 1$ and moreover to the case $1 < p=q < \infty$. In this paper, to address this gap, we provide a proof specifically for $p = q = 1$ using Fubini's theorem and for the case $1 < p=q < \infty $ using different method. We prove that $T_\rho^\alpha$ is bounded from $L^{1, \phi}(\R, d\mu_\alpha)$ to $L^{1, \psi}(\R, d\mu_\alpha)$, while $\tilde{T}_\rho^\alpha$ is bounded from $L^{1, \phi}(\R, d\mu_\alpha)$ to $\mathcal{L}^{1, \psi}(\R, d\mu_\alpha)$, under some suitable conditions on $\rho$, $\phi$ and $\psi$. We also prove that, for $1 < p < \infty$, $T_\rho^\alpha$ is bounded from $L^{p, \phi}(\R, d\mu_\alpha)$ to $L^{p, \psi}(\R, d\mu_\alpha)$, while $\tilde{T}_\rho^\alpha$ is bounded from $L^{p, \phi}(\R, d\mu_\alpha)$ to $\mathcal{L}^{p, \psi}(\R, d\mu_\alpha)$, further we show $\tilde{T}_\rho^\alpha$ is bounded from $\mathcal{L}^{p,\phi}(\R, d\mu_\alpha)$ to $\mathcal{L}^{p, \psi}(\R, d\mu_\alpha)$, under proper conditions on $\rho$, $\phi$ and $\psi$. 

Now the paper is organized as follows. Section \ref{sec:prelim} provides a brief review of Dunkl theory on the real line along with some known results. Section \ref{sec:3} introduces the generalized Dunkl-type Morrey space and the generalized Dunkl-type Campanato space and main results. Section \ref{sec:4} and Section \ref{sec:5} are dedicated to proving the theorems of Section \ref{sec:3}.

\section{\bf{Preliminaries}} \label{sec:prelim}
 The Dunkl operator associated with  a fixed real number $\alpha \geq -\frac{1}{2}$ is defined by
 \begin{align*}
     \Lambda_\alpha f(x) = \frac{df}{dx}(x) + \frac{2 \alpha +1}{x} \frac{f(x) - f(-x)}{2},  ~x \in \mathbb{R},
 \end{align*}
 where $f \in C^\infty(\mathbb{R})$. When $\alpha = -\frac{1}{2}$ the Dunkl derivative $\Lambda_{-\frac{1}{2}}$ is equal to the classical derivative $\frac{d}{dx}$.\\
 For $\alpha \geq - 1/2$ and $\lambda \in \mathbb{C}$, the initial value problem $\Lambda_\alpha f(x) = \lambda f(x)$; $f(0) =1$ has a unique solution $E_\alpha(\lambda x)$ called Dunkl kernel \cite{dunkl1991integral, de1993dunkl, sifi2002generalized}, which is  given by 
 \begin{align*}
     E_\alpha(\lambda x) = \mathcal{J}_\alpha (\lambda x) + \frac{\lambda x}{2(\alpha + 1)} \mathcal{J}_{\alpha+1} (\lambda x), ~ x \in \mathbb{R}, 
 \end{align*}
 where
 \begin{align*}
     \mathcal{J}_\alpha (\lambda x) := \Gamma (\alpha +1) \sum_{n=0}^\infty \frac{(\lambda x / 2)^{2n}}{n! \Gamma(n + \alpha + 1)}
 \end{align*}
 is the modified spherical Bessel function of order $\alpha$. Note that $E_{- 1/ 2}(\lambda x) = e^{\lambda x}$. \\
 The weighted Lebesgue measure $\mu_\alpha$ on $\mathbb{R}$ is given by 
 $$
 d\mu_\alpha(x) :=  A_\alpha |x|^{2 \alpha + 1} dx, \hspace{1cm}  where ~~  A_\alpha = (2 ^ {\alpha + 1 } \Gamma(\alpha + 1))^{-1}.
 $$
 For $1 \leq p \leq \infty$, the space $L^p(\mathbb{R}, d\mu_\alpha)$ denotes the class of  complex-valued measurable functions $f$ on $\mathbb{R}$ such that $$ ||f||_{L^p(\mathbb{R}, d\mu_\alpha)} :=\left(\int_{\mathbb{R}} |f(x)|^p d\mu_\alpha(x)\right)^{\frac{1}{p}} < \infty ,~ ~ if ~ ~ p \in [1, \infty)$$
 and $$ ||f||_{L^{\infty}(\mathbb{R}, d\mu_\alpha)} := ess \displaystyle \sup_{x \in \mathbb{R}} |f(x)| < \infty, ~ ~ if ~ ~ p = \infty.$$\\
Using the Dunkl kernel, one can define the Dunkl transform on $\mathbb{R}$ (see \cite{de1993dunkl}) as follows:\\
 The Dunkl transform of a function $f \in L^1(\mathbb{R}, d\mu_\alpha)$, is given by 
 \begin{align*}
     \mathcal{F}_\alpha (f)(\lambda) := \int_{\mathbb{R}} E_\alpha(- i \lambda x) f(x) d\mu_\alpha(x), ~ \lambda \in \mathbb{R}.
 \end{align*}
 The above integral makes sense as $|E_\lambda (i x)| \leq 1$, $\forall x \in \mathbb{R}$ \cite{rosler1994bessel}. Note that $\mathcal{F}_{- 1/ 2}$ agrees with definition of the classical Fourier transform $\mathcal{F}$, which is given by
 \begin{align*}
     \mathcal{F}(f)(\lambda) := (2 \pi)^{-1/2} \int_{\mathbb{R}} e^{-i \lambda x} f(x) dx, ~ ~ \lambda \in \mathbb{R}.
 \end{align*}
In the below, we list some of the properties of the Dunkl transform $\mathcal{F}_\alpha$. \begin{theorem}
    
\cite{trime2002paley}
 \begin{itemize}
     \item[(i)] For all $ f \in L^1(\mathbb{R}, d\mu_\alpha)$, one has $|| \mathcal{F}_\alpha (f)||_{L^\infty(\mathbb{R}, d\mu_\alpha)} \leq ||f||_{L^1(\mathbb{R}, d\mu_\alpha)}$.
     \item [(ii)] For all $f \in \mathcal{S}(\R)$,
     \begin{align*}
         \mathcal{F}_\alpha ( \Lambda_\alpha f)(\lambda) = i \lambda \mathcal{F}_\alpha (f)(\lambda), ~~ \lambda \in \mathbb{R}.
     \end{align*}
     \item [(iii)] The Dunkl transform $\mathcal{F}_\alpha$ is an isometric isomorphism from $L^2(\mathbb{R}, d\mu_\alpha)$ onto $L^2(\mathbb{R}, d\mu_\alpha)$. In particular, it satisfies the Plancherel formula i.e.,
     \begin{align*}
         ||\mathcal{F}_\alpha (f)||_{L^2(\mathbb{R}, d\mu_\alpha)} = ||f||_{L^2(\mathbb{R}, d\mu_\alpha)}.
     \end{align*}
     \item[(iv)] If $f \in L^1(\mathbb{R}, d\mu_\alpha)$ with $\mathcal{F}_\alpha(f) \in L^1(\mathbb{R}, d\mu_\alpha)$, then one has the inversion formula 
     \begin{align*}
         f(x) = \int_{\mathbb{R}} \mathcal{F}_\alpha (f) (\lambda) E_\alpha(i \lambda x) d\mu_\alpha(\lambda), ~~ a.e.~~ x \in \mathbb{R}.
     \end{align*}
 \end{itemize}
 \end{theorem} 
Now we will denote $B(x,R) = \{y \in \mathbb{R} : |y| \in \,]max\{0, |x|-R\}, |x| + R[\,\}$, $R > 0 $, if $x\neq 0$ and $B(0, R) = \,]-R, R[\, $. Then $\mu_\alpha (B(0, R)) = b_\alpha R^{d_\alpha}$, where $b_\alpha = [2^{\alpha + 1}(\alpha + 1) \Gamma(\alpha + 1)]^{-1}$ and $d_\alpha = 2 \alpha + 2$. We observe that when $|x| \leq r$ then $B(x,r) = B(0, |x| + r) \subseteq B(0, 2r)$ and when $|x| > r$ then $B(x,r)= \{ y \in \mathbb{R} : |y| \in \,]|x|-r, |x|+r[\, \}\subseteq B(0, |x| +r) \subseteq B(0, 2|x|)$.\\
Next for $x, y, z \in \mathbb{R}$, we put
 \begin{align}
 W_\alpha(x, y, z) = (1 - \sigma_{x, y, z} + \sigma_{z, x, y} + \sigma_{z, y, x} ) \kappa_\alpha(|x|, |y|, |z|)  \label{bessel}
 \end{align}
 where 
\begin{align}
 \sigma_{x, y, z} &:=
 \begin{cases}
     \frac{x^2 + y^2 - z^2}{2xy}, & if ~ x, y  \in \mathbb{R}\setminus \{0\},\\
     0, & otherwise, 
 \end{cases}\nonumber
\end{align} 
    and $\kappa_\alpha$ is the Bessel kernel given by 
\begin{align}    
    \kappa_\alpha(|x|, |y|, |z|) :=
    \begin{cases}
        D_\alpha \frac{([(|x| + |y|)^2 - z^2][z^2 - (|x| - |y|)^2])^{\alpha - \frac{1}{2}}}{|xyz|^{2\alpha}}, & if ~ |z| \in S_{x, y},\\
        0, & otherwise,
    \end{cases} \nonumber
\end{align}
Where $D_\alpha = (\Gamma(\alpha + 1))^2/ (2^{\alpha - 1} \sqrt{\pi} \Gamma(\alpha + \frac{1}{2}))$ and $ S_{x, y} = [||x|-|y||,|x|+|y|]$.
Furthermore, if $|z| \in S_{x,y}$ then $|\sigma_{x,y,z}| \leq 1$,$|\sigma_{z,x,y}| \leq 1$ and $|\sigma_{z,y,x}| \leq 1$ ( see \cite{rosler1994bessel}).
\begin{remark} \cite{Rachdi} 
\begin{itemize}
    \item[(i)] The Bessel kernel $\kappa_\alpha$ is symmetric in the three variables, that is, if $x_1, x_2, x_3 \in ]0,+\infty[$ and $\sigma$ is a permutation of $\{1, 2, 3\}$, then
    \begin{align}
        \kappa_\alpha(x_{\sigma(1)}, x_{\sigma(2)}, x_{\sigma(3)}) = \kappa_\alpha(x_1, x_2, x_3). \nonumber
    \end{align}
\item[(ii)] For all $x, y \in \R$,
    \begin{align}
        \int_{0}^{+\infty} \kappa_\alpha(|x|, |y|, z)  d\mu_\alpha(z) = 1. \label{2.1(ii)}
    \end{align}
\end{itemize}
\end{remark}
Using Bessel kernel $\kappa_\alpha$, for $x, y \in \R$ and $f$ a sufficiently smooth function on $[0,\infty[$, one can define the generalized translation $\T_{|x|}^\alpha$ on Bessel–Kingman
hypergroups \cite{Rachdi} given by
\begin{align}
      \T_{|x|}^\alpha f(|y|) = \int_{0}^\infty f(z) \kappa_\alpha(|x|,|y|, z) d\mu_\alpha(z) = \int_{||x|-|y||}^{|x|+|y|} f(z) \kappa_\alpha(|x|,|y|, z) d\mu_\alpha(z).   \label{2.1ac}
\end{align}
 By using $\T_{|x|}^\alpha$, the convolution of two functions $f$ and $g$ on Bessel–Kingman hypergroups is defined by
\begin{align}
      f*_\alpha^h g(|x|) := \int_{0}^\infty f(y) \T_{|x|}^\alpha g(y) d\mu_\alpha(y), \hspace{1cm} x \in \R,  \nonumber
\end{align}
where $*_\alpha^h$ denotes the convolution of two function on the Bessel-Kingman hypergroups associated with $\alpha$. We are using this notation to distinguish it from definition of the Dunkl convolution (see \eqref{eq:dc}). The convolution is associative and commutative.  The translation operator  $\T^\alpha_{|x|}$ satisfies the following properties:
\begin{prop} \cite{Rachdi} \label{prop:1}
\begin{enumerate}
    \item[(i)] For every $x \in \R$, $\T^\alpha_{|x|}$ is a positive-bounded operator on $L^p(\R,d\mu_\alpha)$, $p \in [1, +\infty]$ and for every $f \in L^p(\R, d\mu_\alpha)$,
    \begin{equation}
        \| \T^\alpha_{|x|}f \|_{L^p(\R, d\mu_\alpha)} \leq \| f \|_{L^p(\R,d\mu_\alpha)}. \nonumber
    \end{equation}
    \item[(ii)]  Support of $\T_{|x|}^\alpha \chi_{J(0,r)}$ is contained in the  $J(|x|,r)$  for all $x \in \R$,
where we denote by $J(|x|, r)$ the following set:
\begin{align*}
    J(|x|, r) = ]\max\{0, |x| - r\}, |x| + r [. \nonumber
\end{align*}
\end{enumerate}
\end{prop}
Also note that it follows trivially from \eqref{2.1(ii)}
 and \eqref{2.1ac} that $\T_{|x|}^\alpha \chi_{J(0,r)}(y) \leq 1$.

Next using $W_\alpha(x,y,z)$, the signed measure $\nu_{x, y}$ on $\mathbb{R}$ is given by
  \begin{align}
      d\nu_{x, y}(z) &= \begin{cases}
      W_\alpha(x, y, z) d\mu_\alpha(z), & if ~ x, y \in \R \backslash \{0\},\\ d\delta_x(z), & if ~ y = 0,\\ d\delta_y(z), & if ~ x = 0,
  \end{cases} \nonumber
\end{align} 
where the signed kernel  $W_\alpha$ is an even function with respect to all variables satisfying the following properties ( see \cite{rosler1994bessel}):$$ W_\alpha(x, y, z) = W_\alpha(y, x, z) = W_\alpha(-x, z, y),$$
  $$ W_\alpha(x, y, z) = W_\alpha(-z, y, -x) = W_\alpha(-x, -y, -z)$$ and $$ \int_{\mathbb{R}} |W_\alpha(x, y, z)| d\mu_\alpha(z) \leq 4.$$
     Let $x, y \in \R$ and $f$ be a continuous function on $\R$. Then the Dunkl translation operator of $f$ is defined by
     \begin{align}
         \tau_x^\alpha f(y) = \int_{\R} f(z) d\nu_{x,y}(z).  \label{2.1a}
     \end{align}
    %    Let $z= (x, y)_\theta = \sqrt{x^2 + y^2 - 2|xy| cos \theta}$. Then by change of variable we also have (see \cite{abdelkefi2007characterization} ) 
    % \begin{align}
    %     \tau_x^\alpha f(y) := \frac{c_\alpha}{2} \int_0^\pi &\{[f((x,y)_\theta) + f(-(x,y)_ \theta)]  \nonumber
    %     \\
    %     &+\frac{x+y}{(x,y)_\theta} [f((x,y)_\theta) - f(-(x,y)_\theta)]\}(1 - cos \theta) sin^{2 \alpha}\theta d\theta, \label{2.1b}
    % \end{align}
    % where $ c_\alpha =  \left(\int_0^\pi sin^{2 \alpha}\theta d\theta \right)^{-1} =\frac{\Gamma(\alpha + 1)}{\sqrt{\pi}\Gamma(\alpha + \frac{1}{2})}. $\\
     The Dunkl translation operator satisfies the following properties:
      \begin{prop} \cite{mourou2001transmutation} \label{pro: 2.1}
        \begin{itemize}
            \item [(i)] $\tau_x^\alpha, ~x \in \mathbb{R}$ is a bounded linear operator on $C^\infty(\mathbb{R})$
            \item [(ii)] For all $f \in C^\infty(\mathbb{R})$ and $x,y \in \mathbb{R}$, one has 
            \begin{align*}
                \tau_x^\alpha f(y) = \tau_y^\alpha f(x), \tau_0^\alpha f(x) = f(x), ~ \tau_x^\alpha \circ \tau_y^\alpha = \tau_y^\alpha \circ \tau_x^\alpha.
            \end{align*} 
        \end{itemize}
      \end{prop}

 If $f, g \in L^2(\mathbb{R}, d\mu_\alpha)$, then 
            \begin{align}
                \int_{\mathbb{R}} \tau_x^\alpha f(y) g(y) d\mu_\alpha(y) = \int_{\mathbb{R}} f(y) \tau_{-x}^\alpha g(y) d\mu_\alpha(y), ~ \forall x \in \mathbb{R}. \nonumber
            \end{align}
     
    %  and for $f = f_e + f_o$, where $f_e$ and $f_o$ represents the even and odd part of $f$ respectively, $\tau_x^\alpha$ can be expressed in the following form ( see \cite{rosler1994bessel})
    %  \begin{align}
    %      \tau_x^\alpha f(y) = c_\alpha \int_0^\pi f_e(\sqrt{x^2 + y^2 - 2|xy| cos \theta}) h_1(x, y, \theta) sin^{2\alpha} \theta d\theta \nonumber\\ + c_\alpha \int_0^\pi f_o(\sqrt{x^2 + y^2 - 2|xy| cos \theta}) h_2(x, y, \theta) sin^{2\alpha} \theta d\theta, 
    %      \nonumber
    %   \end{align}
    %  where $$
    %        c_\alpha = \left(\int_0^\pi sin^{2 \alpha}\theta d\theta \right)^{-1} = \frac{\Gamma(\alpha + 1)}{\sqrt{\pi}\Gamma(\alpha + \frac{1}{2})},
    %        $$
    %        \begin{align}
    %            h_1(x, y,\theta) = 1 - sgn(xy) cos\theta  \nonumber
    %        \end{align}
    % and 
    % \begin{align}
    %     h_2(x, y, \theta) = \begin{cases}
    %         \frac{(x+y)[1 - sgn (xy) cos\theta]}{\sqrt{x^2 + y^2 - 2|xy| cos \theta}}, & if ~xy \neq 0\\
    %         0, & if ~xy = 0.
    %     \end{cases}\nonumber
    %     \end{align}
    %     If $f$ is even function on $\mathbb{R}$, then $\tau_x^\alpha$ can be written as
    %     \begin{align*}
    %         \tau_x^\alpha f(y) = c_\alpha \int_{0}^\pi & f(\sqrt{x^2 +y^2 - 2|xy| cos\theta}) sin^{2\alpha} \theta d\theta  
    %         \\
    %          &- sgn(xy). c_\alpha \int_{0}^\pi f(\sqrt{x^2 + y^2 - 2|xy| cos \theta}) cos \theta sin^{2\alpha} \theta d\theta.
    %     \end{align*}
  
Suppose $f$, $g$ are continuous functions on $\mathbb{R}$ with compact support. The generalized convolution $f*_\alpha g$ is defined by 
\begin{align}
    f*_\alpha g(x) := \int_{\mathbb{R}} f(y) \tau_y^\alpha g(x) d\mu_\alpha(y). \label{eq:dc}
\end{align}
The generalized convolution $*_\alpha$ is associative and commutative \cite{rosler1994bessel}. We also mention in the following some useful properties of the generalized convolution $*_\alpha$.
\begin{theorem}\cite{soltani2004lp} \label{tm: 2.1}
\begin{enumerate}
    \item[(i)]   For all $x \in \mathbb{R}$ and $p \geq 1$, the generalized translation operator $\tau_x^\alpha$ is bounded from $L^p(\mathbb{R}, d\mu_\alpha)$ to $L^p(\mathbb{R}, d\mu_\alpha)$ i.e.
      \begin{align*}
          ||\tau_x^\alpha f||_{L^p(\mathbb{R}, d\mu_\alpha)} \leq 4 ||f||_{L^p(\mathbb{R}, d\mu_\alpha)}.
       \end{align*}
    \item [(ii)] If $f\in L^{1}(\mathbb{R}, d\mu_\nu)$, then 
\begin{eqnarray*}
\mathcal{F}_{\alpha} (\tau_x^\alpha f)(\lambda)=E_{\alpha}(i\lambda x)\mathcal{F}_{\alpha}(f)(\lambda),~x,\lambda\in\mathbb{R}.
\end{eqnarray*}
\item [(iii)] If  $f\in L^{1}(\mathbb{R}, d\mu_\alpha)$ and $g\in L^{2}(\mathbb{R}, d\mu_\alpha)$, then 
\begin{eqnarray*}
\mathcal{F}_{\alpha}(f\ast_\alpha g)=\mathcal {F}_{\alpha} (f)\mathcal{F}_{\alpha} (g).
\end{eqnarray*}
 \item[(iv)] Assume that $p,q,r \in [1, \infty]$ satisfying $\frac{1}{p} + \frac{1}{q} = 1 + \frac{1}{r}$  (the Young condition). Then the map $(f,g) \to f *_\alpha g$ defined on $C_c(\mathbb{R}) \times C_c(\mathbb{R})$, extends to a continuous map from
 $L^p(\mathbb{R}, d\mu_\alpha) \times L^q(\mathbb{R}, d\mu_\alpha)$ to $L^r(\mathbb{R}, d\mu_\alpha)$ and one has
\begin{align*}
    ||f *_\alpha g ||_{L^r(\mathbb{R}, d\mu_\alpha)} \leq 4 ||f||_{L^p(\mathbb{R}, d\mu_\alpha)} ||g||_{L^q(\mathbb{R}, d\mu_\alpha)}.
\end{align*} 
\end{enumerate}
\end{theorem}

Now we will recall a few definitions and the known facts in the following, which will be useful in our paper.
\begin{definition} \cite{sps}
  For $1\leq p \leq q < \infty$, the Dunkl-type Morrey space $L^{p,q}(\R, d\mu_\alpha)$ is defined to be the class of all locally $p$-integrable functions $f$ on $\mathbb{R}$ such that \begin{align}
    ||f||_{L^{p,q}(\R, d\mu_\alpha)} :=  \sup_{\substack{{r>0}\\{x \in \R}}} r^{d_\alpha(\frac{1}{q}- \frac{1}{p})} \left(\int_{B(0,r)} \tau_x^\alpha |f|^p(y) d\mu_\alpha(y)\right)^{\frac{1}{p}} < \infty. \label{eq: 3.1a}
\end{align}  
\end{definition}
\begin{definition} \cite{spsa}
    For a function $\phi : (0, +\infty) \to (0, + \infty)$, we define $BMO_{\phi}(\mathbb{R}, d\mu_\alpha)$ to be the space of all locally integrable functions $f$ on $\mathbb{R}$ such that 
\begin{align}
    ||f||_{BMO_{\phi}(\mathbb{R}, d\mu_\alpha)} := \sup_{\substack{{r>0}\\{x \in \R}}} \frac{1}{\phi(r)} \frac{1}{\mu_{\alpha}(B(0,r))} \int_{B(0,r)} |\tau_x^\alpha f(y) -f^\alpha_{B(0,r)}(x) | d\mu_\alpha(y) < \infty, \label{eq: 7.1}
\end{align}
where 
\begin{align}
f^\alpha_{B(0,r)}(x) = \frac{1}{\mu_\alpha(B(0,r))} \int_{B(0,r)} \tau_x^\alpha f(y) d\mu_\alpha(y). \nonumber
\end{align}
If $\phi(r)\equiv 1$, then $BMO_\phi(\mathbb{R}, d\mu_\alpha)= BMO(\mathbb{R}, d\mu_\alpha)$( see \cite{guliyev2009fractional}). 
\end{definition}

     \begin{lemma} \cite{lk}\label{7.1}
      If $f \in L^1(\mathbb{R}, d\mu_\alpha)$ and $g \in L^p(\mathbb{R}, d\mu_\alpha)$, $1 \leq p < \infty$,
then
\begin{align}
\tau_{x_0}^\alpha(f *_\alpha g) = \tau_{x_0}^\alpha f *_\alpha g = f *_\alpha \tau_{x_0}^\alpha g,  ~~x_0 \in \mathbb{R}. \nonumber
\end{align}
    \end{lemma}
    \begin{lemma} (\cite{abdelkefi2007dunkl}, \cite{guliyev2010p}) \label{lem: 2.2}
        Support of $\tau_x^\alpha \chi_{B(0,r)}$ is contained in the ball $B(x,r)$  for all $x \in \R$. Moreover, 
        \begin{align}
            0 \leq \tau_x^\alpha \chi_{B(0, r)}(y) \leq min \left\{ 1, \frac{2 c_\alpha}{2\alpha + 1} \left( \frac{r}{|x|} \right)^{2 \alpha + 1} \right \} , ~\forall~ y \in B(x, r), \nonumber
        \end{align}
        where $c_\alpha = \frac{ \Gamma(\alpha + 1)}{\sqrt{\pi} \Gamma(\alpha + \frac{1}{2})}$. In fact, when $|x| \leq r$, one has $0 \leq \tau_x^\alpha \chi_{B(0,r)}(y) \leq 1$ and when $|x| > r$, one has $0 \leq \tau_x^\alpha \chi_{B(0,r)}(y) \leq \frac{2 c_\alpha}{2 \alpha +1} (\frac{r}{|x|})^{2 \alpha + 1}$.
    \end{lemma}
    % {Also for the even functions $f(x)=f(-x) \in \mathcal{S}(\R)$, the following specific formula of the Dunkl translation was given by R\"{o}sler (\cite{rsl3}).
    % \begin{align}
    %     \tau_x^\alpha f(y) = \int_\R f(\sqrt{x^2 + y^2 - 2x\eta}) d\mu_y^\alpha(\eta) \label{eq: 2.555}
    % \end{align}
    %  where for each $y$ in $\mathbb{R}$, $d\mu_y^\alpha$ is a probability measure on $\mathbb{R}$, whose support is contained in $[-y, y]$.
      \begin{lemma} \cite{spsa} \label{lem: l7.1}
        For $2|x| \leq |y|$, the following inequalities are valid:
        \begin{itemize}
         \item [(i)] $
            \left|\tau_x^\alpha \frac{\rho(|y|)}{|y|^{d_\alpha}} - \frac{\rho(|y|)}{|y|^{d_\alpha}} \right| \leq C |x| \frac{\rho(|y|)}{|y|^{d_\alpha + 1}}, 
            $\\
            \item [(ii)] $
            |\tau_x^\alpha \frac{\rho(|y|)}{|y|^{d_\alpha}}| \leq C \frac{\rho(|y|)}{|y|^{d_\alpha}}, 
             $
         \end{itemize}
        where $\rho$ satisfies the doubling condition  $\frac{1}{2} \leq \frac{r}{s} \leq 2 \implies \frac{1}{C} \leq \frac{\phi(r)}{\phi(s)} \leq C$ and  $C > 0$ is constant.
    \end{lemma}
Now for locally integrable functions $f$ on $\mathbb{R}$, the Dunkl-type Hardy-Littlewood maximal operator $M^\alpha$ ( see \cite{ guliyev2009fractional} ) is defined by
\begin{align}
    M^\alpha f(x) := \sup_{r > 0} \frac{1}{\mu_\alpha(B(0,r))} \int_{B(0,r)} \tau_x^\alpha|f|(y) d\mu_\alpha(y), ~ 
   x \in \mathbb{R}.  \label{eq: 3.4a}
\end{align}
\begin{theorem} \cite{abdelkefi2007dunkl, solatani2005, guliyev2009fractional} \label{tm: 2.2}
    \begin{itemize}
        \item[(i)] Let $s > 0$ be a real number. Then for all $ f \in L^1(\mathbb{R}, d\mu_\alpha) $, one has 
        \begin{align*}
            \mu_\alpha\{ x \in \mathbb{R} : M^\alpha f(x) > s \} \leq \frac{A_1}{s} \int_{\mathbb{R}} |f(x)| d\mu_\alpha(x),
        \end{align*}
       where $A_1 > 0$ is independent of $f$.
    \item[(ii)]  If $f \in L^p(\mathbb{R}, d\mu_\alpha), 1 < p \leq \infty $, then $M^\alpha f \in L^p(\mathbb{R}, d\mu_\alpha)$ and 
    \begin{align*}
        ||M^\alpha f ||_{L^p(\mathbb{R}, d\mu_\alpha)} \leq A_2 ||f||_{L^p(\mathbb{R}, d\mu_\alpha)},
    \end{align*}
    where $A_2 > 0$ is independent of $f$.
    \end{itemize}
\end{theorem}

 \begin{theorem} \cite{spsa} \label{tm: 3.1ba}
For $1 < p < \infty$, let $\phi$ be a positive measurable function on $\mathbb{R}^+$. 
Then for $ f \in L^{p, \phi}(\R, d\mu_\alpha)$, we have 
\begin{align}
    ||M^\alpha f||_{L^{p, \phi}(\R, d\mu_\alpha)} \leq C ||f||_{L^{p, \phi}(\R, d\mu_\alpha)}.\label{eq: 3.6a}
\end{align}
\end{theorem}

\section{\bf{ Generalized fractional integral operators on generalized Morrey spaces, Campanato spaces and main results}} \label{sec:3}
For a function $\rho :(0, +\infty) \to (0, + \infty)$,
\begin{align}
    T_{\rho}^\alpha f(x) := \int_\R  f(y) \tau_x^\alpha \frac{\rho(|y|)}{|y|^{d_\alpha}} d\mu_\alpha(y), \label{eq:3.1}
\end{align}
it is the generalized Dunkl-type fractional integral operator $T_{\rho}^\alpha$. We consider the following conditions on $\rho:$
\begin{align}
      \int_0^1 \frac{\rho(t)}{t} dt < \infty,  \label{eq:3.2}\\
    \frac{1}{2} \leq \frac{r}{s} \leq 2 \implies \frac{1}{C_1} \leq \frac{\rho(r)}{\rho(s)} \leq C_1, \label{eq:3.3}
\end{align} 
where $C_1$ is independent of $r,s > 0$. We define the modified version of the generalized Dunkl-type fractional integral operator $T_\rho^\alpha$ by
\begin{align}
    \tilde{T}_\rho^\alpha f(x) : = \int_{\mathbb{R}} f(y) \left( \tau_x^\alpha \frac{\rho(|y|)}{|y|^{d_\alpha}} - \frac{\rho(|y|)}{|y|^{d_\alpha}} (1- \chi_{B_0}(y)) \right) d\mu_\alpha(y), \label{eq:3.4}
\end{align}
where $B_0 = B(0,1)$ denotes the ball center at the origin and of radius $1$ and  $\chi_{B_0}$ is the characteristic function of $B_0 = B (0,1)$. In this definition, we assume that $\rho$ satisfies  \eqref{eq:3.2}, \eqref{eq:3.3} and the following conditions: 
    \begin{align}
         \int_r^\infty \frac{\rho(t)}{t^2} dt \leq C'_1 \frac{\rho(r)}{r} ~~ for ~~ all~~ r > 0, \label{eq:3.5}\\
       \frac{1}{2} \leq \frac{r}{s} \leq 2 \implies \left|\frac{\rho(r)}{r^{d_\alpha}} - \frac{\rho(s)}{s^{d_\alpha}}\right| \leq C''_1 |r- s| \frac{\rho(s)}{s^{d_\alpha + 1}}, \label{eq:3.6}
    \end{align}
    where $C'_1, C''_2 > 0$ are independent of $r,s > 0$.
    For example, the function $\rho(r) = r^\beta$ satisfies \eqref{eq:3.2}, \eqref{eq:3.3} and \eqref{eq:3.6} for $0 < \beta < d_\alpha$, and also satisfies \eqref{eq:3.5} for $ 0 < \beta < 1. $
Here $K_\rho^\alpha(x) = \frac{\rho(|x|)}{|x|^{d_\alpha}}$ is known as the generalized Dunkl-type Riesz kernel. If $\rho(t) = t^{\beta}$, $0< \beta < d_\alpha$, then we have the Dunkl type fractional integral operator $I_{\beta}^\alpha = T_{\rho}^\alpha$ with Dunkl-type Riesz kernel $K_{\rho}^\alpha = K_{\beta}^\alpha$.
  \begin{remark}
     Since the modified version of the generalized fractional integral operator $T_\rho^\alpha $ \eqref{eq:3.4} can be written as 
    \begin{align}
        \tilde{T}_\rho^\alpha f(x) = f *_\alpha \frac{\rho(|.|)}{|.|^{d_\alpha}}(x) - f *_\alpha \frac{\rho(|.|)}{|.|^{d_\alpha}}(1 - 
  \chi_{B_0}) (0), \nonumber 
    \end{align}
    hence by using Lemma \ref{7.1}, we can write 
    \begin{align}
        \tau_{x_0}^\alpha \tilde{T}_\rho^\alpha f(x) &= \tau_{x_0}^\alpha f *_\alpha \frac{\rho(|.|)}{|.|^{d_\alpha}} (x) - \tau_{x_0}^\alpha f *_\alpha \frac{\rho(|.|)}{|.|^{d_\alpha}}(1 - \chi_{B_0}) (0) \nonumber
        \\
        &= \int_{\mathbb{R}} \tau_{x_0}^\alpha f(y) \left( \tau_x^\alpha \frac{\rho(|y|)}{|y|^{d_\alpha}} - \frac{\rho(|y|)}{|y|^{d_\alpha}}(1 - \chi_{B_0}(y)) \right) d\mu_\alpha(y). \label{eq: 7a}
    \end{align} 
     \end{remark}
      For functions $ \theta_1, \theta_2 :(0, +\infty) \to (0, + \infty)$, we denote $\theta_1(r) \sim \theta_2(r)$ if there exists a constant $C>0$ such that $ \frac{1}{C}\theta_1(r) \leq \theta_2(r) \leq C \theta_1(r)$, $r>0$.
      
 A function $\theta : (0, +\infty) \to (0, + \infty)$ is said to be almost increasing (almost decreasing) if there exists a constant $C>0$ such that $\theta(r) \leq C \theta(s) ~ (\theta(r) \geq C\theta(s))$  for $r \leq s$.
 
Let $\phi : \R^{+} \to \R^{+}$ be a function which is measurable. For $1\leq p < \infty$, we define the generalized Dunkl-type Morrey space $L^{p, \phi}(\R, d\mu_\alpha)$ for all locally $p$-integrable functions $f$ on $\mathbb{R}$ such that 
\begin{align}
    ||f||_{L^{p, \phi}(\R, d\mu_\alpha)} :=  \sup_{\substack{{r > 0}\\{x \in \R}}} \frac{1}{\phi(r)}\left(\frac{1}{\mu_\alpha(B(0,r))}\int_{B(0,r)} \tau_x^\alpha |f|^p(y) d\mu_\alpha(y) \right)^{\frac{1}{p}} < \infty. \label{eq:3.7}
\end{align} 

For $1 \leq p < \infty$ and a function $\phi : (0, +\infty) \to (0, + \infty)$, we define the generalized Dunkl-type Campanato space $\mathcal{L}^{p,\phi}(\mathbb{R}, d\mu_\alpha)$ to be the space of all locally integrable functions $f$ on $\mathbb{R}$ such that 
\begin{align}
    ||f||_{\mathcal{L}^{p,\phi}(\mathbb{R}, d\mu_\alpha)} := \sup_{\substack{{r>0}\\{x \in \R}}} \frac{1}{\phi(r)} \left(\frac{1}{\mu_{\alpha}(B(0,r))} \int_{B(0,r)} |\tau_x^\alpha f(y) -f^\alpha_{B(0,r)}(x) |^p d\mu_\alpha(y) \right)^{\frac{1}{p}} < \infty, \label{eq:3.8}
\end{align}
where 
\begin{align}
f^\alpha_{B(0,r)}(x) = \frac{1}{\mu_\alpha(B(0,r))} \int_{B(0,r)} \tau_x^\alpha f(y) d\mu_\alpha(y). \nonumber
\end{align}

 In order to define $L^{p, \phi}(\R, d\mu_\alpha)$ and $\mathcal{L}^{p,\phi}(\R, d\mu_\alpha)$, we assume that $\phi$ satisfies the doubling condition \eqref{eq:3.3}and $r^{d_\alpha} \phi^p(r)$ is almost increasing. For $p=1$, $\mathcal{L}^{1,\phi}(\R,d\mu_\alpha)= BMO_{\phi}(\R, d\mu_\alpha)$, which is defined in \eqref{eq: 7.1}.
 One can show that $ f $ belongs to $\mathcal{L}^{p,\phi}(\R, d\mu_\alpha)$ if there exists a constant $C < \infty $ such that, for each  $ r > 0$ and $x \in \R $, a constant $c_{B(0,r)}(x) < \infty $ satisfying  
\begin{align*}
 \frac{1}{\phi(r)} \left( \frac{1}{\mu_\alpha(B(0,r))} \int_{B(0,r)} |\tau_x^\alpha f(y) - c_{B(0,r)}(x)|^p d\mu_\alpha(y) \right)^{1/p} < C.
\end{align*}
This implies that $\|f\|_{\mathcal{L}^{p,\phi}}$ is bounded by $ 2C $. We will use this result in the proof of our results.

In the following, we state the main theorems of our paper.
% If $\phi(r)\equiv 1$, then $BMO_\phi(\mathbb{R}, d\mu_\alpha)= BMO(\mathbb{R}, d\mu_\alpha)$( see \cite{guliyev2009fractional}). 
\begin{theorem} \label{3.1}
    Let $\phi$ and $\psi$ satisfy the doubling condition \eqref{eq:3.3}. Let $\rho$ satisfy \eqref{eq:3.2} and \eqref{eq:3.3}. Assume that there exists a constant $ C > 0$ such that, for all $r > 0$, 
    \begin{align}
        \int_0^r \frac{\rho(t)}{t} dt ~\phi(r) + \int_r^\infty \frac{\rho(t) \phi(t)}{t} dt \leq C \psi(r) \label{eq:3.9}.
    \end{align}
    Then $T_\rho^\alpha$ is bounded from $L^{1, \phi}(\R, d\mu_\alpha)$ to $L^{1, \psi}(\R, d\mu_\alpha)$.
\end{theorem}
\begin{theorem} \label{3.3}
    Let $\phi$ and $\psi$ satisfy the doubling condition \eqref{eq:3.3}. Let $\rho$ satisfy \eqref{eq:3.2} and \eqref{eq:3.3}. Assume that there exists a constant $ C > 0$ such that, for all $r > 0$, 
    \begin{align}
        \int_0^r \frac{\rho(t)}{t} dt ~\phi(r) + \int_r^\infty \frac{\rho(t) \phi(t)}{t} dt \leq C \psi(r) \label{eq:3.11}.
    \end{align}
    Then for each $1 < p < \infty$, $T_\rho^\alpha$ is bounded from $L^{p, \phi}(\R, d\mu_\alpha)$ to $L^{p, \psi}(\R, d\mu_\alpha)$.
\end{theorem}
\begin{theorem} \label{3.2}
   Let $\phi$ and $\psi$ satisfy the doubling condition \eqref{eq:3.3}. Let $\rho$ satisfy \eqref{eq:3.2}, \eqref{eq:3.3}, \eqref{eq:3.5} and \eqref{eq:3.6}. Assume that there exists a constant $ C > 0$ such that, for all $r > 0$, 
    \begin{align}
        \int_0^r \frac{\rho(t)}{t} dt ~\phi(r) + r\int_r^\infty \frac{\rho(t) \phi(t)}{t^2} dt \leq C \psi(r) \label{eq:3.10}.
    \end{align}
    Then $\tilde{T}_\rho^\alpha$ is bounded from $L^{1, \phi}(\R, d\mu_\alpha)$ to $\mathcal{L}^{1, \psi}(\R, d\mu_\alpha)$.  
\end{theorem}
 % \begin{theorem}\label{eq: t8}
 %      Let $\phi$ and $\psi$ satisfy the doubling condition \eqref{eq:3.3}. Let $\rho$ satisfy \eqref{eq:3.2}, \eqref{eq:3.3}, \eqref{eq:3.5} and \eqref{eq:3.6}. Assume that there exists a constant $ C > 0$ such that, for all $r > 0$, 
 %    \begin{align}
 %        \int_0^r \frac{\rho(t)}{t} dt ~\phi(r) + r\int_r^\infty \frac{\rho(t) \phi(t)}{t^2} dt \leq C \psi(r) \label{eq:3.10}.
 %    \end{align}
 %    Then $\tilde{T}_\rho^\alpha$ is bounded from $\mathcal{L}^{1,\phi}(\mathbb{R}, d\mu_\alpha)$ to $\mathcal{L}^{1,\psi}(\mathbb{R}, d\mu_\alpha)$.
 %  \end{theorem}
 %  \begin{Cor} (\cite{spsa})
 %  Let $\phi$ and $\psi$ satisfy the doubling condition \eqref{eq:3.3} and be almost increasing. Let $\rho$ satisfy \eqref{eq:3.2}, \eqref{eq:3.3}, \eqref{eq:3.5} and \eqref{eq:3.6}. Assume that there exists a constant $ C > 0$ such that, for all $r > 0$, 
 %    \begin{align}
 %        \int_0^r \frac{\rho(t)}{t} dt ~\phi(r) + r\int_r^\infty \frac{\rho(t) \phi(t)}{t^2} dt \leq C \psi(r) \label{eq:3.10}.
 %    \end{align}
 %   Then $\tilde{T}_\rho^\alpha$ is bounded from $BMO_\phi(\mathbb{R}, d\mu_\alpha)$ to $BMO_\psi(\mathbb{R}, d\mu_\alpha)$.
 %  \end{Cor}

\begin{theorem} \label{3.4}
   Let $\phi$ and $\psi$ satisfy the doubling condition \eqref{eq:3.3}. Let $\rho$ satisfy \eqref{eq:3.2}, \eqref{eq:3.3} and \eqref{eq:3.6}. Assume that there exists a constant $ C > 0$ such that, for all $r > 0$, 
    \begin{align}
        \int_0^r \frac{\rho(t)}{t} dt ~\phi(r) + r\int_r^\infty \frac{\rho(t) \phi(t)}{t^2} dt \leq C \psi(r) \label{eq:3.12}.
    \end{align}
    Then for each $1 < p < \infty$, $\tilde{T}_\rho^\alpha$ is bounded from $L^{p, \phi}(\R, d\mu_\alpha)$ to $\mathcal{L}^{p, \psi}(\R, d\mu_\alpha)$.  
\end{theorem}
\begin{theorem} \label{3.5}
   Let $\phi$ and $\psi$ satisfy the doubling condition \eqref{eq:3.3} and $\int_1^\infty \frac{\phi(t)}{t} dt < \infty$. Let $\rho$ satisfy \eqref{eq:3.2}, \eqref{eq:3.3}, \eqref{eq:3.5} and \eqref{eq:3.6}. Assume that there exists a constant $ C > 0$ such that, for all $r > 0$, 
    \begin{align}
       \int_r^\infty \frac{\phi(t)}{t} dt \int_0^r \frac{\rho(t)}{t} dt  + r\int_r^\infty \frac{\rho(t) \phi(t)}{t^2} dt \leq C \psi(r) \label{eq:3.13}.
    \end{align}
    Then for each $1 < p < \infty$, $\tilde{T}_\rho^\alpha$ is bounded from $\mathcal{L}^{p, \phi}(\R, d\mu_\alpha)$ to $\mathcal{L}^{p, \psi}(\R, d\mu_\alpha)$.  
\end{theorem}
%%%%%%%%%%%%%%%%%%%%%%%%%%%%%%%%%%%%%%%%%%%%%%%%%%%%%%%%%%%%%%%%%%%%%%%%%%%%%%%%%%%%%%%%%%%%%%%%%%%%%%%%%%%%%%%%%%%%%%%%%%%%%%%%%%%%%%%%%%%%%%%%%%%%%
\section{\bf{Proofs of the boundedness of generalized fractional integral operators on generalized Morrey spaces}} \label{sec:4}
In this section, we provide the proofs of Theorem \ref{3.1} and Theorem \ref{3.3}. Towards this, we first prove the following proposition and lemma. We remark that Proposition \ref{prop:4.1b} establishes a relation between Dunkl translation $\tau_x^\alpha$ and generalized translation on the Bessel–Kingman hypergroups $\T_{|x|}^\alpha$, which we will be using often in proving our results.

 \begin{prop} \label{prop:4.1b}
      For all $x,y \in \mathbb{R}$ and for any non-negative function $f \in L^p(\R,d\mu_\alpha)$, $p \in [1, +\infty)$,
    \begin{align}
        | \tau_x^\alpha f(y) |^{p} \leq 8^p\T_{|x|}^\alpha f_e^p(|y|). \label{eq:23}
    \end{align}
     \end{prop}
     
\textbf{Proof of Proposition \ref{prop:4.1b}.}   From \eqref{bessel}, we have
  \begin{align*}
|W_\alpha(x,y,z)| \leq 4 \kappa_\alpha(|x|,|y|,|z|), \hspace{1cm}  |z| \in S_{x,y}.
  \end{align*}
  Now decomposing $f$ as $f = f_e + f_o$ with $f_e$ even part and $f_o$ odd part, we get from \eqref{2.1a}
    \begin{align}
        | \tau_x^\alpha f (y)| &\leq 4 \int_{S_{x,y} \cup -S_{x,y}} f(z) \kappa_\alpha(|x|,|y|,|z|) d\mu_\alpha(z) \nonumber\\
        &\leq 8 \int_{||x| - |y||}^{|x| + |y|} f_e(z) \kappa_\alpha(|x|, |y|, |z|) \, d\mu_\alpha(z) \nonumber \\
        & = 8\int_0^\infty f_e(z) \kappa_\alpha(|x|, |y|, |z|) \, d\mu_\alpha(z) = 8\T_{|x|}^\alpha f_e(|y|), \nonumber
    \end{align}
    where the last equality follows from the definition of the generalized translation $\T_{|x|}^\alpha$ on Bessel-Kingman hypergroups defined in \eqref{2.1ac}. Therefore, result is clear for $p=1$. Now 
    assume that $p \in (1, +\infty)$ and let $q$ be the conjugate exponent of $p$, i.e., $\frac{1}{p} + \frac{1}{q} = 1$. We write
    \begin{align}
|f_e(z) \kappa_\alpha(|x|, |y|, |z|)| = f_e(z) \kappa_\alpha(|x|, |y|, |z|)^{1/p}\kappa_\alpha(|x|, |y|, |z|)^{1/q}.
    \end{align}
    Applying H\"older's inequality and using \eqref{2.1(ii)} , we obtain
    \begin{align}
        |\tau_x^\alpha f(y) |^p \leq 8^p \int_0^\infty f_e^p(z) \kappa_\alpha(|x|,|y|,|z|) d\mu_\alpha(z) = 8^p\T_{|x|}^\alpha f_e^p(|y|). \nonumber 
    \end{align} \qed

     \begin{prop} \label{prop:4.2}
       Let $f \in L^1(\R,d\mu_\alpha)$ be any function. Then
   \begin{align}
    \int_{I(y, \rho)} f(x) d\mu_\alpha(x) =\int_{B(0,\rho)} \tau_{y}^\alpha f(x) d\mu_\alpha(x), \nonumber
\end{align}
for all $y \in \R$ and for all $\rho > 0$.
     \end{prop}
     
\textbf{Proof of Proposition \ref{prop:4.2}.} For $I(y, \rho)= (y -\rho, y+ \rho)$, consider 
\begin{align}
     \int_{I(y, \rho)} f(x) d\mu_\alpha(x) &=  \int_{\R} f(x) \chi_{I(y, \rho)}(x)  d\mu_\alpha(x) \nonumber
     \\
     & = \int_{\R} f(x) \chi_{I(-x, \rho)}(-y)  d\mu_\alpha(x) \nonumber
     \\
      & = \int_{\R} f(x) \tau_0^\alpha \chi_{I(-x, \rho)}(-y)  d\mu_\alpha(x).  \label{eq: 4.3ab}
     \end{align}
     where we have used the fact $\chi_I(y,\rho)(x)= \chi_{I(x, \rho)}(y)= \chi_{I(-x, \rho)}(-y)$.  Using Proposition \ref{pro: 2.1} in \eqref{eq: 4.3ab}, we have
     \begin{align}
      \int_{I(y, \rho)} f(x) d\mu_\alpha(x) & = \int_{\R} f(x) \tau_{-y}^\alpha \chi_{I(-x, \rho)}(0)  d\mu_\alpha(x) \nonumber
     \\ 
      & = \int_{\R} \tau_y^\alpha f(x) \chi_{I(-x, \rho)}(0)  d\mu_\alpha(x) \nonumber
     \\
      & = \int_{\R} \tau_y^\alpha f(x) \chi_{I(0, \rho)}(x)  d\mu_\alpha(x) \nonumber
     \\
       & = \int_{B(0, \rho)} \tau_y^\alpha f(x)  d\mu_\alpha(x), \nonumber
\end{align}
where again we have used the fact $\chi_{I(-x, \rho)}(0)= \chi_{I(0, \rho)}(x)$ and $I(0,\rho)= B(0,\rho)$.  \qed
\begin{lemma} \label{lem:4.1}
    Let $\rho, \phi$ and $\psi$ be functions from $\R^+$ to $\R^+$ with the doubling condition and $\phi(r) r^{d_\alpha}$ be almost decreasing. If there exists a constant $ C > 0$ such that for all $r > 0$ 
    \begin{align}
        r^j \int_r^{\infty} \frac{\rho(t) \phi(t)}{t^{1 + j}} dt \leq C \psi(r), \label{eq:4.1}
    \end{align}
    where $j=0$ or $j=1$, then there exist a constant $C' > 0$ such that, for all $f \in L^{1, \phi}(\R, d\mu_\alpha)$, $ x_0, x \in R$ and $r > 0$
    \begin{align}
        r^j \int_{B^c(0,r)} |\tau_{x_0}^\alpha f(y)| \tau_x^\alpha \frac{\rho(|y|)}{|y|^{d_\alpha + j}} d\mu_\alpha(y) \leq C' \psi(r) ||f||_{L^{1, \phi}(\R, d\mu_\alpha)}. \label{eq: 4.2}
    \end{align}
    % In addition, for all $f \in \mathcal{L}^{1, \phi}(\R, d\mu_\alpha)$, we also have
    % \begin{align}
    %     r^j\int_{B^c(0,r)} |\tau^\alpha_{x_0} f(y) - f_{B(0,r)}(x_0)| \tau^\alpha_x \frac{\rho(y)}{|y|^{d_\alpha +j}} d\mu_\alpha(y) \leq C'' \psi(r) ||f||_{\mathcal{L}^{1, \phi}(\R, d\mu_\alpha)}.
    % \end{align}
\end{lemma}

\textbf{Proof of Lemma \ref{lem:4.1}.} For $j=0$ or $j=1$, consider
\begin{align}
    & r^j  \int_{B^c(0, r)} |\tau_{x_0}^\alpha f(y)| \tau_x^\alpha \frac{\rho(|y|)}{|y|^{d_\alpha + j}}  d\mu_\alpha(y) \nonumber\\
    & = r^j \sum_{k = 0}^{\infty} \int_{2^k r \leq |y| < 2^{k+1}r} |\tau_{x_0}^\alpha f(y)| \tau_x^\alpha\frac{\rho(|y|)}{|y|^{d_\alpha + j}}   d\mu_\alpha(y) \nonumber
   \\
     & \leq C r^j \sum_{k= 0}^{\infty}   \tau_x^\alpha\frac{\rho(2^k r)}{( 2^k r)^{d_\alpha + j}}\int_{2^k r \leq |y| < 2^{k+1}r} |\tau_{x_0}^\alpha f(y)|  d\mu_\alpha(y) \nonumber
     \\
      & = C r^j \sum_{k= 0}^{\infty}  \frac{\rho(2^k r)}{( 2^k r)^{d_\alpha + j}}\int_{2^k r \leq |y| < 2^{k+1}r} |\tau_{x_0}^\alpha f(y)|  d\mu_\alpha(y). \label{eq:26b}
    \end{align}
    For each $k = 0, 1, 2, ...,$ we will first estimate 
\begin{align}
    \int_{2^k r \leq |y| < 2^{k+1} r} |\tau_{x_0}^\alpha f(y)| d\mu_\alpha(y). \nonumber
\end{align}
Since the operator $\tau_{x_0}^\alpha$ is linear it suffices to estimate the above integral for non-negative functions only. Then Proposition \ref{prop:1} and \eqref{eq:23} lead to
\begin{align}
    \int_{2^k r \leq |y| < 2^{k+1} r} |\tau_{x_0}^\alpha f(y)| d\mu_\alpha(y) &\leq 8 \int_{0}^{2^{k+1}r} \T_{|x_0|}^\alpha f_e(|y|)  d\mu_\alpha(y)  \nonumber
    \\
    & = 8 \int_0^\infty f_e(y) \T_{|x_0|}^\alpha \chi_{(0, 2^{k+1}r)}(y) d\mu_\alpha(y) \nonumber
    \\
    & = 8 \int_{J(|x_0|, 2^{k+1} r)} f_e(y) \T_{|x_0|}^\alpha \chi_{J(0, 2^{k+1}r)}(y) d\mu_\alpha(y) \nonumber
    \\
      & \leq 8 \int_{J(|x_0|, 2^{k+1} r)} f_e(y) d\mu_\alpha(y), \label{eq:26cc}
\end{align}
where $f = f_e +f_0$ with $f_e$ even and $f_0$ odd.\\
We observe that 
\begin{align}
J(|x_0|, 2^{k+1}r) \subseteq I(|x_0|, 2^{k+1}r)= \left(|x_0| -2^{k+1}r, |x_0| + 2^{k+1}r\right). \nonumber
\end{align}
Hence from \eqref{eq:26cc}, we obtain
\begin{align}
     \int_{2^k r \leq |y| < 2^{k+1} r} |\tau_{x_0}^\alpha f(y)| d\mu_\alpha(y) \leq C \int_{I(|x_0|, 2^{k+1}r)} f_e(y) d\mu(y). \label{eq:26c}
\end{align}
Now using  Proposition \ref{prop:4.2}), 
% \begin{align}
%     \int_{I(y, \rho)} |f|(x) d\mu_\alpha(x) \leq C\int_{B(0,\rho)} \tau_{y}^\alpha |f|(x) d\mu_\alpha(x), \nonumber
% \end{align}
% for all $y \in \R$ and for all $\rho > 0$,
we get from \eqref{eq:26c}
% the above step can be interpreted using the concept of the space of homogeneous type. This means that $X = \R$ is a topological space equipped with a continuous quasi metric $\rho$ and a positive measure $\mu$ satisfying the doubling property
% \begin{align}
% \mu(E(x, 2r)) \leq C_0 \mu(E(x,r)) \nonumber
% \end{align}
% for a constant $C_0$ independent on $x \in X$ and $r >0$, where $E(x,r) = \{ y \in X : \rho(x,y) < r \}, ~ \rho(x,y)= |x-y|$ is open ball in $X$. Note that the Dunkl measure $\mu_\alpha$  satisfies the doubling condition and hence $(\R, |x-y|, d\mu_\alpha)$ is a space of homogeneous type. We obtain from \eqref{eq:26c}
% \begin{align}
%      \int_{2^k r \leq |y| < 2^{k+1} r} |\tau_{x_0}^\alpha f(y)| d\mu_\alpha(y) \leq C \int_{E(x_0, 2^{k+1}r)} f_e(y) d\mu(y). \nonumber
% \end{align}
% Since, by using the result of Amri and Sifi (\cite{amri2012riesz}), i.e. $\min\{|x_0-y|, |x_0+y|\} \leq \sqrt{x_0^2 + y^2 -2x_0\eta} \leq \max\{|x_0-y|, |x_0+y|\}$, $\forall~ x_0, x \in \R$, $\forall~ \eta \in [-y,y]$ in \eqref{eq: 2.555}, it follows that 
% \begin{align}
%     \int_{E(x_0, 2^{k+1}r)} f_e(y) d\mu(y) \leq C\int_{B(0,2^{k+1}r)} \tau_{x_0}^\alpha f_e(y) d\mu_\alpha(y). \nonumber
% \end{align}
% Using the above fact, we get
\begin{align}
    \int_{2^k r \leq |y| < 2^{k+1} r} |\tau_{x_0}^\alpha f(y)| d\mu_\alpha(y) & \leq C \int_{B(0, 2^{k+1} r)} \tau_{|x_0|}^\alpha f_e(y) d\mu_\alpha(y) \nonumber \\
    & \leq C \phi(2^k r) (2^k r)^{d_\alpha} ||f_e||_{L^{1,\phi}(\mathbb{R}, d\mu_\alpha)} \nonumber
    \\
    & \leq C \phi(2^k r) (2^k r)^{d_\alpha} ||f||_{L^{1,\phi}(\mathbb{R}, d\mu_\alpha)}. 
     \label{eq: 4.5d}
\end{align}
Thus \eqref{eq:26b} reduces to
\begin{align}
     r^j  \int_{B^c(0, r)} |\tau_{x_0}^\alpha f(y)| \tau_x^\alpha \frac{\rho(|y|)}{|y|^{d_\alpha + j}}  d\mu_\alpha(y) & \leq C r^j ||f||_{L^{1, \phi}(\mathbb{R}, d\mu_\alpha)} \sum_{k= 0}^{\infty} \frac{\rho(2^k r)}{( 2^k r)^{d_\alpha + j}} \phi(2^{k} r) (2^k r)^{d_\alpha} \nonumber
     \\
     & \leq C r^j ||f||_{L^{1, \phi}(\mathbb{R}, d\mu_\alpha)} \sum_{k= 0}^{\infty} \frac{\rho(2^{k+1} r) \phi(2^{k+1} r)}{(2^k r)^{j }} \nonumber
     \\
     &\leq C r^j ||f||_{L^{1, \phi}(\R, d\mu_\alpha)} \sum_{k=0}^{\infty} \int_{2^k r}^{2^{k+1} r} \frac{\rho(t) \phi(t)}{t^{j + 1}} dt \nonumber
      \\
     &\leq C r^j ||f||_{L^{1, \phi}(\R, d\mu_\alpha)} \int_{r}^{\infty} \frac{\rho(t) \phi(t)}{t^{j + 1}} dt, \nonumber
\end{align}
since
\begin{align}
    \int_{2^k r}^{2^{k+1}r} \frac{\rho(t)\phi(t)}{t^{j + 1}} dt &\geq C\frac{\rho(2^{k+1}r)\phi(2^{k+1}r)}{(2^{k+1}r)^{j +1}} 2^k r \geq C\frac{\rho(2^{k+1}r)\phi(2^{k+1} r)}{(2^k r)^{\gamma - d_\alpha }} . \nonumber
\end{align}
Finally using \eqref{eq:4.1}, we obtain
\begin{align}
    r^j  \int_{B^c(0, r)} |\tau_{x_0}^\alpha f(y)| \tau_x^\alpha \frac{\rho(|y|)}{|y|^{d_\alpha + j}}  d\mu_\alpha(y) \leq C' \psi(r) ||f||_{L^{1, \phi}(\R, d\mu_\alpha)}.  \nonumber 
\end{align}
This completes the proof. \qed
% The proof of the second estimate follows similarly as in \cite{spsa}, the proof of Lemma 7.2], by replacing $||f||_{BMO_\phi(\R,d\mu_\alpha)}$ with $||f||_{\mathcal{L}^{1,\phi}(\R,d\mu_\alpha)}$. \qed

\textbf{Proof of Theorem \ref{3.1}.}
Using the definition of generalized Dunkl-type fractional integral operator $T_\rho^\alpha$ (see \eqref{eq:3.1}) and Lemma \ref{7.1}, for any $x_0 \in \R$ and $r>0$ we can write
\begin{align}
     \int_{B(0,r)} \tau_{x_0}^\alpha |T_{\rho}^\alpha f|(x)  d\mu_\alpha(x) 
     &= \int_{\R}  |T_\rho^\alpha f|(x) \tau_{-x_0}^\alpha \chi_{B(0,r)}(x)d\mu_\alpha(x). \nonumber
     \\
     &\leq \int_{\R} \left(\int_{\R} |f|(y) \tau_x^\alpha \frac{\rho(y)}{|y|^{d_\alpha}} d\mu_\alpha(y) \right) \tau_{-x_0}^\alpha \chi_{B(0,r)}(x)d\mu_\alpha(x) \nonumber
     \\
    &= \left| \int_{\R} \left(\int_{\R} |f|(y) \tau_x^\alpha \frac{\rho(y)}{|y|^{d_\alpha}} d\mu_\alpha(y) \right) \tau_{-x_0}^\alpha \chi_{B(0,r)}(x)d\mu_\alpha(x) \right| \nonumber
     \\
     &= \left|\int_{B(0,r)} \left( \int_{\R} \tau_{x_0}^\alpha |f|(y) \tau_x^\alpha \frac{\rho(|y|)}{|y|^{d_\alpha}} d\mu_\alpha(y) \right) d\mu_\alpha(x) \right| \nonumber
     \\
     &\leq \int_{B(0,r)} \left( \int_{\R} |\tau_{x_0}^\alpha |f|(y)| \tau_x^\alpha \frac{\rho(|y|)}{|y|^{d_\alpha}} d\mu_\alpha(y) \right) d\mu_\alpha(x)  \nonumber
     \\
     &= \int_{B(0,r)} T_1(x) d\mu_\alpha(x) + \int_{B(0,r)} T_2(x) d\mu_\alpha(x), \label{eq:8.1}
\end{align}
where 
\begin{align}
T_1(x) := \int_{B(0,r)}  |\tau_{x_0}^\alpha |f|(y)| \tau_x^\alpha \frac{\rho(|y|)}{|y|^{d_\alpha}} d\mu_\alpha(y) \nonumber
\end{align}
and
\begin{align}
    T_2(x) :=\int_{B^c(0,r)} |\tau_{x_0}^\alpha |f|(y)| \tau_x^\alpha \frac{\rho(|y|)}{|y|^{d_\alpha}} d\mu_\alpha(y). \nonumber 
 \end{align}
Now consider
\begin{align}
    &\int_{B(0, r)}\left( \int_{B(0,r)} |\tau_{x_0}^\alpha |f|(y)| \tau_x^\alpha \frac{\rho(|y|)}{|y|^{d_\alpha}} d\mu_\alpha(x) \right ) d\mu_\alpha(y) \nonumber\\
    & = \int_{B(0, r)} |\tau_{x_0}^\alpha |f|(y)| \left( \int_{B(0, r)} \tau_y^\alpha \frac{\rho(|x|)}{|x|^{d_\alpha}} d\mu_\alpha(x) \right) d\mu_\alpha(y). \nonumber 
    \end{align}
    Using the support property of the Dunkl translation, it follows at once that
    \begin{align}
     &\int_{B(0, r)}\left( \int_{B(0,r)} |\tau_{x_0}^\alpha |f|(y)| \tau_x^\alpha \frac{\rho(|y|)}{|y|^{d_\alpha}} d\mu_\alpha(x) \right ) d\mu_\alpha(y) \nonumber\\
     & = \int_{B(0, r)} |\tau_{x_0}^\alpha |f|(y)| \left( \int_{B(-y, r)} \frac{\rho(|x|)}{|x|^{d_\alpha}} \tau_{-y}^\alpha \chi_{B(0,r)}(x) d\mu_\alpha(x) \right) d\mu_\alpha(y) \nonumber \\
     & \leq \int_{B(0, r)} |\tau_{x_0}^\alpha |f|(y)| \left( \int_{B(0, 2r)} \frac{\rho(|x|)}{|x|^{d_\alpha}} d\mu_\alpha(x) \right) d\mu_\alpha(y) \nonumber \\
    & = \int_{B(0, r)} |\tau_{x_0}^\alpha |f|(y)| \left( \int_0^{2r} \frac{\rho(t)}{t} dt \right) d\mu_\alpha(y). \nonumber 
    \end{align}
    Then Proposition \ref{prop:1} and Proposition \ref{prop:4.1b} with $p=1$ lead to
    \begin{align}
     &\int_{B(0, r)}\left( \int_{B(0,r)} |\tau_{x_0}^\alpha |f|(y)| \tau_x^\alpha \frac{\rho(|y|)}{|y|^{d_\alpha}} d\mu_\alpha(x) \right ) d\mu_\alpha(y) \nonumber\\
     & \leq 16 \int_0^r \T^\alpha_{|x_0|} |f|_e(|y|) \left( \int_0^{2r} \frac{\rho(t)}{t} dt \right) d\mu_\alpha(y)  \nonumber \\
     & = 16\int_0^\infty |f|_e(y)  \T^\alpha_{|x_0|} \chi_{(0,r)}(y)  \left( \int_0^{2r} \frac{\rho(t)}{t} dt \right) d\mu_\alpha(y) \nonumber\\
     &= 16\int_{J(|x_0|, r)} |f|_e(y)  \T^\alpha_{|x_0|} \chi_{J(0,r)}(y)  \left( \int_0^{2r} \frac{\rho(t)}{t} dt \right) d\mu_\alpha(y)\nonumber\\
     &\leq 16\int_{J(|x_0|, r)} |f|_e(y)  \left( \int_0^{2r} \frac{\rho(t)}{t} dt \right) d\mu_\alpha(y) , \label{eq:34c}
    \end{align}
    where $|f|_e$ denotes the even part of $|f|$.\\
    Now, in the above step, using the same observation  as one introduced in the proof of Lemma \ref{lem:4.1}, we have
    \begin{align}
        \int_{J(|x_0|,r)} |f|_e(y)| d\mu_\alpha(y) \leq C \int_{I(|x_0|,r)} |f|_e(y) d\mu(y) = C \int_{B(0,r)} \tau_{|x_0|}^\alpha |f|_e(y) d\mu_\alpha(y). \label{eq:34d}
    \end{align}
    Using the above \eqref{eq:34d} in \eqref{eq:34c}, we get
\begin{align}
&\int_{B(0, r)}\left( \int_{B(0,r)} |\tau_{x_0}^\alpha |f|(y)| \tau_x^\alpha \frac{\rho(|y|)}{|y|^{d_\alpha}} d\mu_\alpha(x) \right ) d\mu_\alpha(y) \nonumber
\\
&\leq C \int_{B(0, r)} \tau_{|x_0|}^\alpha |f|_e(y)  \left( \int_0^{2r} \frac{\rho(t)}{t} dt \right) d\mu_\alpha(y). \label{eq:34e}
\end{align}
 Given that $\rho$ satisfies the doubling condition, using the definition of generalized Dunkl-type Morrey space (see\eqref{eq:3.7}) and using \eqref{eq:3.9}, we obtain from \eqref{eq:34e}
    \begin{align}
    & \int_{B(0, r)}\left( \int_{B(0,r)} |\tau_{x_0}^\alpha |f|(y)| \tau_x^\alpha \frac{\rho(|y|)}{|y|^{d_\alpha}} d\mu_\alpha(x) \right ) d\mu_\alpha(y) \nonumber \\
     &\leq C \phi(r) r^{d_\alpha} ||f_e||_{L^{1,\phi}(\R, d\mu_\alpha)} \int_0^{2r} \frac{\rho(t)}{t} dt \nonumber\\
      &\leq C ||f_e||_{L^{1,\phi}(\mathbb{R}, d\mu_\alpha)} r^{d_\alpha} \phi(r) \int_0^r \frac{\rho(t)}{t} dt.\nonumber\\
     &\leq C ||f||_{L^{1,\phi}(\mathbb{R}, d\mu_\alpha)} r^{d_\alpha} \psi(r).  \label{34a} 
\end{align}
Hence Fubini's theorem leads to
\begin{align}
    \int_{B(0, r)} T_1(x) d\mu_\alpha(x) \leq C \psi(r) r^{d_\alpha} ||f||_{L^{1,\phi}(\mathbb{R}, d\mu_\alpha)}. \label{eq:4.3}
\end{align}
Now since
\begin{align}
    T_2(x) = \int_{B^c(0,r)} |\tau_{x_0}^\alpha |f|(y)|  \tau_x^\alpha \frac{\rho(|y|)}{|y|^{d_\alpha}} d\mu_\alpha(y) , \nonumber
\end{align}
therefore Lemma \ref{lem:4.1} with $j=0$ implies  
\begin{align}
  |T_2(x)| \leq C \psi(r) ||f||_{L^{1,\phi}(\mathbb{R} , d\mu_\alpha)}. \nonumber
\end{align}
This gives
\begin{align}
    \int_{B(0,r)} T_2(x) d\mu_\alpha(x) \leq C \psi(r) r^{d_\alpha}||f||_{L^{1,\phi}(\mathbb{R} , d\mu_\alpha)}.  \label{eq:4.4}
\end{align}
Finally, by substituting  the values from \eqref{eq:4.3} and \eqref{eq:4.4} in \eqref{eq:8.1}, we have 
\begin{align}
    \int_{B(0,r)} \tau_{x_0}^\alpha |T_\rho^\alpha f|(x) d\mu_\alpha(x) \leq C ||f||_{L^{1,\phi}(\mathbb{R}, d\mu_\alpha)} r^{d_\alpha} \psi(r), \nonumber
\end{align} 
which implies
\begin{align}
    ||T_\rho^\alpha f ||_{L^{1,\psi}(\mathbb{R}, d\mu_\alpha)}   \leq C ||f||_{L^{1,\phi}(\mathbb{R}, d\mu_\alpha)}, \nonumber
\end{align}
completing the proof.
\qed

\begin{remark}
    Note that the authors proved the boundedness of $T_\rho^\alpha$ from  $L^{p, \phi}(\R,d\mu_\alpha)$ to $L^{q,\psi}(\R, d\mu_\alpha)$, for $1 < p<q<\infty$ (see Theorem $6.1$, \cite{spsa}) using the boundedness of maximal operator $M^\alpha$ from $L^{p,\phi}(\R,d\mu_\alpha)$ to $L^{p,\phi}(\R,d\mu_\alpha)$. However, when $p=q=1$, the boundedness of $T_\rho^\alpha$ from $L^{1,\phi}(\R, d\mu_\alpha)$ to $L^{1, \psi}(\R,d\mu_\alpha)$ does not follow from the proof of Theorem $6.1$ in \cite{spsa} since maximal operator is not bounded from $L^{1,\phi}(\R,d\mu_\alpha)$ to $L^{1,\phi}(\R, d\mu_\alpha)$. Therefore, in this article, we provide the proof of $T_\rho^\alpha$ in the case of $p=q=1$, using the generalized translation on the Bessel–Kingman hypergroups.
\end{remark}

\textbf{Proof of Theorem \ref{3.3}}
We know that 
\begin{align}
    T_{\rho}^\alpha f(x) := \int_\R  f(y) \tau_x^\alpha \frac{\rho(|y|)}{|y|^{d_\alpha}} d\mu_\alpha(y). \nonumber
\end{align}
 Let $r > 0$ be given. Then for $f \in {L}^{p, \phi}(\R, d\mu_\alpha)$, we write
 \begin{align}
     |T_{\rho}^\alpha f (x)| &\leq  \left( \int_{B(0,r)} + \int_{B^c(0,r)} \right) |f|(y) \tau_{x}^\alpha \frac{\rho(|y|)}{|y|^{d_\alpha}} d\mu_\alpha(y) \nonumber
   \\
   & =  \left( \int_{B(0,r)} + \int_{B^c(0,r)} \right) \tau_{-x}^\alpha |f|(y) \frac{\rho(|y|)}{|y|^{d_\alpha}} d\mu_\alpha(y) = T_{1}(x) + T_{2}(x) , \nonumber
 \end{align} 
where \begin{align}
    T_{1}(x) := \int_{B(0, r)}  \tau_{-x}^\alpha |f|(y)  \frac{\rho(|y|)}{|y|^{d_\alpha}} d\mu_\alpha(y) \nonumber
    \end{align}
    and
    \begin{align}
    T_{2}(x) := \int_{B^c(0, r)}  \tau_{-x}^\alpha |f|(y)  \frac{\rho(|y|)}{|y|^{d_\alpha}}  d\mu_\alpha(y). \nonumber
\end{align}
 Now, proceeding similarly as in \cite{spsa} and
using \eqref{eq:3.11}, we have
\begin{align}
    T_{1}(x) &\leq C M^\alpha f(-x) \sum_{k=-\infty}^{-1} \int_{2^k r}^{2^{k+1}r} \frac{\rho(t)}{t} dt = C M^\alpha f(-x) \int_0^r \frac{\rho(t)}{t} dt \nonumber
    \\ 
    & \leq C \frac{\psi(r)}{\phi(r)} M^\alpha f(-x) . \label{eq:4.8}
\end{align}
Again, proceeding similarly as in \cite{spsa} and
using \eqref{eq:3.11} for the second term $T_2$, we obtain
\begin{align}
   |T_{2}(x)| &\leq C ||f||_{L^{p,\phi}(\R, d\mu_\alpha)} \int_r^\infty \frac{\rho(t)\phi(t)}{t} dt \nonumber\\
   &\leq C \psi(r) ||f||_{L^{p, \phi}(\R, d\mu_\alpha)}. \label{eq:4.9}
\end{align}
Now for $1 \leq p < \infty$ and for all $r>0$, consider
\begin{align}
    \frac{1}{\psi(r)^p r^{d_\alpha}} \int_{B(0,r)} \tau_{x_0}^\alpha |T_1(x)|^p d\mu_\alpha(x) &= \frac{1}{\psi(r)^p r^{d_\alpha}} \int_{\R} |T_1(x)|^p \tau_{-x_0}^\alpha \chi_{B(0,r)} (x)d\mu_\alpha(x). \label{eq:a}
   \end{align}
   Substituting the estimate for $T_1$ from \eqref{eq:4.8} in the above inequality \eqref{eq:a}, we get
   \begin{align}
   \frac{1}{\psi(r)^p r^{d_\alpha}} \int_{B(0,r)} \tau_{x_0}^\alpha |T_1(x)|^p d\mu_\alpha(x) &\leq \frac{C}{\phi(r)^p r^{d_\alpha}} \int_{\R} (M^\alpha f(x))^p \tau_{-x_0}^\alpha \chi_{B(0,r)}(x) d\mu_\alpha(x). \nonumber 
   \end{align}
    Using the  boundedness of Dunkl-type maximal operator $M^\alpha$ on $L^{p, \phi}(\mathbb{R},d\mu_\alpha)$ from Theorem \ref{tm: 3.1ba}, we get
   \begin{align}
    \frac{1}{\psi(r)^p r^{d_\alpha}} \int_{B(0,r)} \tau_{x_0}^\alpha |T_1(x)|^p d\mu_\alpha(x) &\leq \frac{C}{\phi(r)^p r^{d_\alpha}} \int_{B(0,r)} \tau_{x_0}^\alpha (M^\alpha f(x))^p d\mu_\alpha(x) \nonumber 
    \\
    &\leq  \frac{C}{\phi(r)^p r^{d_\alpha}} \phi(r)^p (r)^{d_\alpha} ||M^\alpha f||_{L^{p, \phi}(\R, d\mu_\alpha)}^p \nonumber\\
    &= C ||f||^p_{L^{p, \phi}(\R, d\mu_\alpha)}
   \label{eq:4.14a}
\end{align}
% since $\phi(r)^p r^{d_\alpha}$ is almost decreasing. Hence, in the both case, we have
% \begin{align}
%     \frac{1}{\psi(r)^p r^{d_\alpha}} \int_{B(0,r)} \tau_{x_0}^\alpha |T_1(x)|^p d\mu_\alpha(x)  \leq C ||f||^p_{L^{p, \phi}(\R, d\mu_\alpha)}. \label{eq:4.14a}
% \end{align}
Again, for the second term,
using \eqref{eq:4.9}, we get
\begin{align}
    \frac{1}{\psi(r)^p r^{d_\alpha}} \int_{B(0,r)} \tau_{x_0}^\alpha|T_2(x)|^p d\mu_\alpha(x) &=  \frac{1}{\psi(r)^p r^{d_\alpha}} \int_{\R} |T_2(x)|^p \tau_{-x_0}^\alpha \chi_{B(0,r)}(x) d\mu_\alpha(x) \nonumber
    \\
    &\leq C ||f||^p_{L^{p,\phi(\R, d\mu_\alpha)}}. \label{eq:4.15a}
\end{align}
Now it is easy to observe that
\begin{align}
    | T_\rho^\alpha f(x)|^p \leq 2^{p-1} (|T_1(x)|^p + |T_2(x)|^p). \nonumber
\end{align}
Therefore combining the two estimates \eqref{eq:4.14a} and \eqref{eq:4.15a}, we obtain
\begin{align}
   & \frac{1}{\psi(r)^p r^{d_\alpha}} \int_{B(0,r)} \tau_{x_0}^\alpha |T_\rho^\alpha f|^p(x) d\mu_\alpha(x) \nonumber \\
   &=  \frac{1}{\psi(r)^p r^{d_\alpha}} \int_{\R} |T_\rho^\alpha f|^p(x)  \tau_{-x_0}^\alpha \chi_{B(0,r)}(x) d\mu_\alpha(x) \nonumber\\
    & \leq  C\frac{1}{\psi(r)^p r^{d_\alpha}} \int_{\R} (|T_1(x)|^p + |T_2(x)|^p)  \tau_{-x_0}^\alpha \chi_{B(0,r)}(x) d\mu_\alpha(x) \nonumber\\
     & =  C\frac{1}{\psi(r)^p r^{d_\alpha}} \int_{B(0,r)} \tau_{x_0}^\alpha |T_1(x)|^p d\mu_\alpha(x) +  \frac{1}{\psi(r)^p r^{d_\alpha}} \int_{B(0,r)} \tau_{x_0}^\alpha |T_2(x)|^p d\mu_\alpha(x) \nonumber\\
   &\leq C ||f||^p_{L^{p,\phi(\R, d\mu_\alpha)}}, \nonumber
\end{align}
and the result follows by taking supremum over $r> 0$ and $x_0 \in \R$ in the definition of generalized Dunkl-type Morrey space.\qed 
%%%%%%%%%%%%%%%%%%%%%%%%%%%%%%%%%%%%%%%%%%%%%%%%%%%%%%%%%%%%%%%%%%%%%%%%%%%%%%%%%%%%%%%%%%%%%%%%%%%%%%%%%%%%%%%%%%%%%%%%%%%%%%%%%%%%%%%%%%%%%%%%%%%%%%%%%%%%%%%%%%%%%%%%%%%%%%%
\section{\bf{Proofs of the boundedness of modified version of generalized fractional integral operators on generalized Campanato spaces}} \label{sec:5}
In this section, we provide the proofs of Theorem \ref{3.2}, Theorem \ref{3.4} and Theorem \ref{3.5}. Towards this, we first introduce the following lemmas. 

\begin{lemma} \label{lem:4.2}
  Let $f \in \mathcal{L}^{p, \phi}(\R, d\mu_\alpha)$, $1 \leq p < \infty$ and $\phi$ satisfies the doubling condition. Then for any ball $B(0,r)$ and for any $x_0 \in \R$ and $k= 1,2,3, ...$, we have 
  \begin{align}
      &\left( \frac{1}{\mu_\alpha(B(0, 2^k r))} \int_{B(0, 2^k r)} |\tau_{x_0}^\alpha f(y) - f^\alpha_{B(0,r)}(x_0)|^p d\mu_\alpha(y) \right)^{1/p} \nonumber
      \\
      &\leq C ||f||_{\mathcal{L}^{p,\phi}(\R, d\mu_\alpha)} \int_r^{2^k r} \frac{\phi(t)}{t} dt, \nonumber
  \end{align}
  where $C > 0$ is dependent only on $d_\alpha$ and the doubling constant of $\phi$. 
\end{lemma}

\begin{lemma} \label{lem:4.3}
   Let $f \in \mathcal{L}^{p, \phi}(\R, d\mu_\alpha)$ for some $1 \leq p < \infty$ and $\phi$ satisfy the doubling condition. If $B(0, r) \subseteq B(0, s)$ in $\R$, then
   \begin{align}
       |f^\alpha_{B(0,r)}(a) - f^\alpha_{B(0,s)}(b)| \leq C ||f||_{\mathcal{L}^{p,\phi}(\R, d\mu_\alpha)} \int_r^{2s} \frac{\phi(t)}{t} dt, \nonumber
   \end{align}
 where $C > 0$ is dependent only on $d_\alpha$ and the doubling constant of $\phi$. 
\end{lemma}

\begin{remark}
      We omit the proof of Lemma \ref{lem:4.2} and Lemma \ref{lem:4.3}, since it follows as in the same line of \cite{eridani2004generalized}.
  \end{remark}
  \begin{lemma}\label{lem:7}
    Let $1 \leq p < \infty$, $\phi$ satisfy the doubling condition and $ \int_1^\infty \frac{\phi(t)}{t} dt < \infty$. If $f \in \mathcal{L}^{p, \phi}(\R, d\mu_\alpha)$, then $f_{B(0,r)}(x_0)$ converges as $r$ tends to infinity and 
    \begin{align}
        ||f - \lim_{r \to \infty} f^\alpha_{B(0,r)}(x_0) ||_{L^{p, \tilde{\phi}}(\R, d\mu_\alpha)} \leq C ||f||_{\mathcal{L}^{p, \phi}(\R, d\mu)}, \nonumber
    \end{align}
    where $\tilde{\phi}(r) = \int_r^\infty \frac{\phi(t)}{t} dt $ and $ C > 0$ is dependent only on $d_\alpha$ and the doubling constant of $\phi$.
\end{lemma}

\textbf{Proof of Lemma \ref{lem:7}}
Since from Lemma \ref{lem:4.3}
\begin{align*}
    |f^\alpha_{B(0,r)}(a) - f^\alpha_{B(0,s)}(b)| \leq C ||f||_{\mathcal{L}^{p,\phi}(\R, d\mu_\alpha)} \int_r^{2s} \frac{\phi(t)}{t} dt,
\end{align*}
it follows that there exists a constant $\sigma(f)$, independent of $a \in \R$,
such that
\begin{align*}
    \lim_{r\to \infty} f^\alpha_{B(0,r)}(a) = \sigma(f),
\end{align*}
and
\begin{align*}
    |f^\alpha_{B(0,r)}(a) - \sigma(f)| \leq C ||f||_{\mathcal{L}^{p,\phi}(\R, d\mu_\alpha)} \int_r^\infty \frac{\phi(t)}{t} dt.
\end{align*}
Hence we have, for all $B=B(0, r)$ and $a \in \R$
\begin{align*}
   & \frac{1}{\mu_\alpha(B(0,r))} \int_{B(0,r)}\tau_a^\alpha |f-\sigma(f)|^p(x) d\mu_\alpha(x) \\
   &= \frac{1}{\mu_\alpha(B(0,r))} \int_{\R} |f- \sigma(f)|^p(x) \tau_{-a}^\alpha \chi_{B(0,r)}(x) d\mu_\alpha(x) \\
   &= \frac{1}{\mu_\alpha(B(0,r))} \int_{\R} |f- f^\alpha_{B(0,r)}(a) + f^\alpha_{B(0,r)}(a) - \sigma (f)|^p \tau_{-a}^\alpha \chi_{B(0,r)}(x) d\mu_\alpha(x) \\
   &\leq \frac{C}{\mu_\alpha(B(0,r))} \Bigg(\int_{\R} |f-f^\alpha_{B(0,r)}(a)|^p \tau_{-a}^\alpha \chi_{B(0,r)}(x) d\mu_\alpha(x) \\
   & \hspace{3cm}+\int_{\R} |f^\alpha_{B(0,r)}(a) - \sigma(f)|^p \tau_{-a} \chi_{B(0,r)}(x) d\mu_\alpha(x) \Bigg).\nonumber
   \end{align*}
Now using the support of $\tau_{-a}^\alpha \chi_{B(0,r)}$ from Lemma \ref{lem: 2.2} and proceeding similarly as in \cite{spsa}, and by the definition of the generalized Campanato norm, we obtain
\begin{align*}
   & \frac{1}{\mu_\alpha(B(0,r))} \int_{B(0,r)}\tau_a^\alpha |f-\sigma(f)|^p(x) d\mu_\alpha(x) \nonumber
   \\
    & \leq  \frac{C}{\mu_\alpha(B(0,r))} \int_{B(-a,r)} |f- f^\alpha_{B(0,r)}(a)|^p d\mu_\alpha(x) + C |f^\alpha_{B(0,r)}(a) - \sigma (f)|^p \\
   & \leq C ||f||^p_{\mathcal{L}^{p,\phi}(\R, d\mu_\alpha)} \phi^p(r) + C ||f||^p_{\mathcal{L}^{p, \phi}} \tilde{\phi}^p(r) \\
   & \leq C ||f||^p_{\mathcal{L}^{p,\phi}(\R,d\mu_\alpha)} \tilde{\phi}^p(r),
\end{align*}
which implies 
\begin{align*}
    ||f -\lim_{r\to \infty} f_{B(0,r)}(a) ||_{L^{p, \tilde{\phi}}(\R, d\mu_\alpha)} \leq C ||f||_{\mathcal{L}^{p, \phi}(\R, d\mu_\alpha)}.
\end{align*} \qed 

\textbf{Proof of Theorem \ref{3.2}.} 
  Let $f \in L^{1,\phi}(\mathbb{R}, d\mu_\alpha)$. For given $r >0$, let $\tilde{B}=B(0,2r)$ and suppose $x \in B(0,r)$.
  Consider
  \begin{align}
     & E_{B(0,r)}(x) = \int_\mathbb{R} \tau_{x_0}^\alpha f(y)  \left( \tau_x^\alpha \frac{\rho(|y|)}{|y|^{d_\alpha}} - \frac{\rho(|y|)}{|y|^{d_\alpha}} (1 - \chi_{\tilde{B}}(y))\right) d\mu_\alpha(y),  \nonumber
     \\
      & C_{B(0,r)}(x_0) = \int_\mathbb{R} \tau_{x_0}^\alpha f(y) \left( \frac{\rho(|y|)}{|y|^{d_\alpha}}(1 - \chi_{\tilde{B}}(y)) - \frac{\rho(|y|)}{|y|^{d_\alpha}} (1 - \chi_{B_0}(y))\right) d\mu_\alpha(y), \nonumber
       \\
         & E_{B^1(0,r)}(x) = \int_{\tilde{B}}  \tau_{x_0}^\alpha f(y)  \tau_x^\alpha \frac{\rho(|y|)}{|y|^{d_\alpha}} d\mu_\alpha(y), \hspace{6.25 cm} \nonumber
         \\
          & E_{B^2(0,r)}(x) = \int_{\tilde{B}^c}  \tau_{x_0}^\alpha f(y)  \left( \tau_x^\alpha \frac{\rho(|y|)}{|y|^{d_\alpha}} - \frac{\rho(|y|)}{|y|^{d_\alpha}}\right) d\mu_\alpha(y),  \nonumber 
  \end{align}
  where $x_0 \in \mathbb{R}$. Then from \eqref{eq: 7a}
  \begin{align}
      \tau_{x_0}^\alpha \tilde{T}_\rho^\alpha f(x) - C_{B(0,r)}(x_0)  = E_{B(0,r)}(x) = E_{B^1(0,r)}(x) + E_{B^2(0,r)}(x). \label{eq:38a}
  \end{align}
  Since 
  \begin{align}
     &\left|\frac{\rho(|y|) (1 - \chi_{\tilde{B}}(y))}{|y|^{d_\alpha}} - \frac{\rho(|y|)(1 - \chi_{B_0}(y))}{|y|^{d_\alpha}}\right| \nonumber\\
     & \hspace{3cm}\leq  \begin{cases}
            0, &  |y| < min (1,2r) ~or ~|y| \geq max(1, 2r);\\
            \frac{\rho(|y|)}{|y|^{d_\alpha}}
            , & otherwise,
        \end{cases}\nonumber
          \end{align}
    this implies $C_{B(0,r)}(x_0)$ is finite.
Now consider
\begin{align}
    &\int_{B(0, 2r)}\left( \int_{B(0,r)} |\tau_{x_0}^\alpha f(y) | \tau_x^\alpha \frac{\rho(|y|)}{|y|^{d_\alpha}} d\mu_\alpha(x) \right ) d\mu_\alpha(y) \nonumber\\
    & \leq \int_{B(0, 2r)} |\tau_{x_0}^\alpha f(y)| \left( \int_{B(0, r)} \tau_y^\alpha \frac{\rho(|x|)}{|x|^{d_\alpha}} d\mu_\alpha(x) \right) d\mu_\alpha(y) \nonumber \\
     & \leq \int_{B(0, 2r)} |\tau_{x_0}^\alpha f(y)| \left( \int_{B(-y, r)} \frac{\rho(|x|)}{|x|^{d_\alpha}} \tau_{-y}^\alpha \chi_{B(0,r)}(x) d\mu_\alpha(x) \right) d\mu_\alpha(y) \nonumber \\
     & \leq \int_{B(0, 2r)} |\tau_{x_0}^\alpha f(y)| \left( \int_{B(0, 3r)} \frac{\rho(|x|)}{|x|^{d_\alpha}} d\mu_\alpha(x) \right) d\mu_\alpha(y) \nonumber \\
    & = \int_{B(0, 2r)} |\tau_{x_0}^\alpha f(y) | \left( \int_0^{3r} \frac{\rho(t)}{t} dt \right) d\mu_\alpha(y). \nonumber
    \end{align}
     Now, proceeding similarly to the proof of Theorem \ref{3.1}, we get from \eqref{34a} 
    \begin{align}
     \int_{B(0, 2r)}&\left(\int_{B(0,r)} |\tau_{x_0}^\alpha f(y) | \tau_x^\alpha \frac{\rho(|y|)}{|y|^{d_\alpha}} d\mu_\alpha(x) \right ) d\mu_\alpha(y) \nonumber
     \\
     &\leq C ||f||_{L^{1,\phi}(\mathbb{R}, d\mu_\alpha)} r^{d_\alpha} \phi(r) \int_0^r \frac{\rho(t)}{t} dt \nonumber\\
     & \leq C||f||_{L^{1,\phi}(\R,d\mu_\alpha)} r^{d_\alpha} \psi(r). \nonumber
\end{align}
Therefore, applying Fubini's theorem, we get 
\begin{align}
    \int_{B(0, r)} |E_{B^1(0,r)}(x)| d\mu_\alpha(x) \leq C \psi(r) r^{d_\alpha} ||f||_{L^{1,\phi}(\mathbb{R}, d\mu_\alpha)}, \nonumber
\end{align}
which in addition shows that $E_{B^1(0,r)}$ is finite a.e.\\
Now we will find an estimate for $E_{B^2(0,r)}$. Towards this, we observe that under hypothesis \eqref{eq:3.10} of Theorem \ref{3.2}, the condition \eqref{eq:4.1} in Lemma \ref{lem:4.1} is satisfied. Therefore, for $x \in B(0,r)$, Lemma \ref{lem: l7.1} and Lemma \ref{lem:4.1} with $j=1$ give
\begin{align}
    |E_{B^2(0,r)}(x)| &\leq \int_{B^c(0,2r)} |\tau_{x_0}^\alpha f(y) | \left| \tau_x^\alpha \frac{\rho(|y|)}{|y|^{d_\alpha}} - \frac{\rho(|y|)}{|y|^{d_\alpha}} \right| d\mu_\alpha(y) \nonumber\\
   &\leq C r \int_{B^c(0,2r)} |\tau_{x_0}^\alpha f(y)| \frac{\rho(|y|)}{|y|^{d_\alpha +1}} d\mu_\alpha(y)  \nonumber\\
  % & \leq C \rho(r) \phi(r) ||f||_{L^{1,\phi}(\mathbb{R}, d\mu_\alpha)} \nonumber\\
  % & \leq C \int_0^r \frac{\rho(t)}{t} dt \phi(r) ||f||_{L^{1,\phi}(\mathbb{R}, d\mu_\alpha)} \nonumber\\
  & \leq C \psi(r) ||f||_{L^{1,\phi}(\mathbb{R} , d\mu_\alpha)}, \nonumber
\end{align}
which also shows that $E_{B^2(0,r)}$ is finite.\\
Hence
\begin{align}
    \int_{B(0,r)} |E_{B^2(0,r)}(x)| d\mu_\alpha(x) \leq C \psi(r) r^{d_\alpha}||f||_{L^{1,\phi}(\mathbb{R} , d\mu_\alpha)}.  \nonumber
\end{align}
Finally, it follows from \eqref{eq:38a} 
\begin{align}
  \frac{1}{\psi(r)} \frac{1}{\mu_\alpha(B(0,r))} \int_{B(0,r)} |\tau_{x_0}^\alpha \tilde{T}_\rho^\alpha f(x) - C_{B(0,r)}(x_0)| d\mu_\alpha(x) \leq C ||f||_{L^{1,\phi}(\mathbb{R}, d\mu_\alpha)}. \nonumber
\end{align} 
Taking supremum over $x_0 \in \R$ and $r >0$, we get from the definition of generalized Dunkl-type Campanato space given in \eqref{eq:3.8}
\begin{align}
    ||\tilde{T}_\rho^\alpha f ||_{\mathcal{L}^{1,\psi}(\mathbb{R}, d\mu_\alpha)}   \leq C ||f||_{L^{1,\phi}(\mathbb{R}, d\mu_\alpha)}, \nonumber
\end{align}
completing the proof.\qed

 % \begin{remark} \label{rem:4.3}
 %         In [the proof of Theorem 7.1, \cite{spsa}] through the Replacement of $f \in BMO_{\phi}(\R, d\mu_\alpha)$ with $f \in \mathcal{L}^{1,\phi}(\R,d\mu_\alpha)$  and by making use of Lemma \ref{lem:4.1}, one can show that $\tilde{T}_\rho^\alpha$ is bounded from $\mathcal{L}^{1,\phi}(\mathbb{R}, d\mu_\alpha)$ to $\mathcal{L}^{1,\psi}(\mathbb{R}, d\mu_\alpha)$.
         
 %         Here, the functions $\phi$ and $\psi$ satisfy the doubling condition \eqref{eq:3.3}, and $\phi(r) r^{d_\alpha}$ and $\psi(r) r^{d_\alpha}$ are almost decreasing. The function $\rho$ satisfies \eqref{eq:3.2}, \eqref{eq:3.3}, \eqref{eq:3.5} and \eqref{eq:3.6}. Furthermore,  assume that there exists a constant $ C > 0$ such that, for all $r > 0$, 
 %    \begin{align}
 %        \int_0^r \frac{\rho(t)}{t} dt ~\phi(r) + r\int_r^\infty \frac{\rho(t) \phi(t)}{t^2} dt \leq C \psi(r) \label{eq:3.10}.
 %    \end{align}
 %    \end{remark}
 
\textbf{Proof of Theorem \ref{3.4}}
 Let $f \in L^{p,\phi}(\mathbb{R}, d\mu_\alpha)$. For given $r >0$, let $\tilde{B}=B(0,2r)$ and suppose $x \in B(0,r)$.
  Consider
  \begin{align}
     & E_{B(0,r)}(x) = \int_\mathbb{R} \tau_{x_0}^\alpha f(y)  \left( \tau_x^\alpha \frac{\rho(|y|)}{|y|^{d_\alpha}} - \frac{\rho(|y|)}{|y|^{d_\alpha}} (1 - \chi_{\tilde{B}}(y))\right) d\mu_\alpha(y),  \nonumber
     \\
      & C_{B(0,r)}(x_0) = \int_\mathbb{R} \tau_{x_0}^\alpha f(y) \left( \frac{\rho(|y|)}{|y|^{d_\alpha}}(1 - \chi_{\tilde{B}}(y)) - \frac{\rho(|y|)}{|y|^{d_\alpha}} (1 - \chi_{B_0}(y))\right) d\mu_\alpha(y), \nonumber
       \\
         & E_{B^1(0,r)}(x) = \int_{\tilde{B}}  \tau_{x_0}^\alpha f(y)  \tau_x^\alpha \frac{\rho(|y|)}{|y|^{d_\alpha}} d\mu_\alpha(y), \hspace{6.25 cm} \nonumber
         \\
          & E_{B^2(0,r)}(x) = \int_{\tilde{B}^c}  \tau_{x_0}^\alpha f(y)  \left( \tau_x^\alpha \frac{\rho(|y|)}{|y|^{d_\alpha}} - \frac{\rho(|y|)}{|y|^{d_\alpha}}\right) d\mu_\alpha(y),  \nonumber 
  \end{align}
  where $x_0 \in \mathbb{R}$. Then from \eqref{eq: 7a}
  \begin{align}
      \tau_{x_0}^\alpha \tilde{T}_\rho^\alpha f(x) - C_{B(0,r)}(x_0)  = E_{B(0,r)}(x) = E_{B^1(0,r)}(x) + E_{B^2(0,r)}(x). \label{eq:51c}
  \end{align}
We have already shown that $C_{B(0,r)}(x_0)$ is finite in the proof of Theorem \ref{3.2}. Now
\begin{align}
    E_{B^1(0,r)}(x) &= \int_{B(0,2r)} \tau_{x_0}^\alpha f(y) \tau_x^\alpha \frac{\rho(|y|)}{|y|^{d_\alpha}} d\mu_\alpha(y) \nonumber\\
    &= \int_{\R} ((\tau_{x_0}^\alpha f)\chi_{B(0,2r)})(y) \tau_x^\alpha \frac{\rho(|y|)}{|y|^{d_\alpha}} d\mu_\alpha(y) \nonumber\\
    &= \int_{B(-x, 2r)} \tau_{-x}^\alpha ((\tau_{x_0}^\alpha f)\chi_{B(0,2r)})(y) \frac{\rho(|y|)}{|y|^{d_\alpha}} d\mu_\alpha(y). \nonumber
\end{align}
We write $\tilde{f} := (\tau_{x_0}^\alpha f)\chi_{B(0,2r)}$  to have
\begin{align}
    E_{B^1(0,r)}(x) = \int_{B(-x, 2r)} \tau_{-x}^\alpha \tilde{f}(y) \frac{\rho(|y|)}{|y|^{d_\alpha}} d\mu_\alpha(y), \nonumber
\end{align}
which implies,
\begin{align}
    |E_{B^1(0,r)}(x)| &\leq \int_{B(-x, 2r)} |\tau_{-x}^\alpha \tilde{f}(y)| \frac{\rho(|y|)}{|y|^{d_\alpha}} d\mu_\alpha(y) \nonumber\\
    & \leq \int_{B(0, 3r)} |\tau_{-x}^\alpha \tilde{f}(y)| \frac{\rho(|y|)}{|y|^{d_\alpha}} d\mu_\alpha(y)  \nonumber \\
    &\leq \int_{B(0, 4r)} |\tau_{-x}^\alpha \tilde{f}(y)| \frac{\rho(|y|)}{|y|^{d_\alpha}} d\mu_\alpha(y) \nonumber \\
     & = \sum_{k= -\infty}^{1} \int_{2^k r \leq |y| < 2^{k+1}r}  |\tau_{-x}^\alpha \tilde{f}(y)|  \frac{\rho(|y|)}{|y|^{d_\alpha}} d\mu_\alpha(y) \nonumber
    \\
    & \leq C \sum_{k= -\infty}^{1} \frac{\rho(2^k r)}{( 2^k r)^{d_\alpha}}\int_{2^k r \leq |y| < 2^{k+1}r}  |\tau_{-x}^\alpha \tilde{f}(y)|  d\mu_\alpha(y). \label{eq: 4.12a}
    \end{align}
Since the operator $\tau_x^\alpha$ is linear, it suffices to estimate the above integral for non-negative functions only. So we assume $\tilde{f}$ to be non-negative. Now from Proposition \ref{prop:4.1b}, we have 
\begin{align}
    \int_{2^k R \leq |y| < 2^{k+1} R} |\tau_{-x}^\alpha \tilde{f}(y)| d\mu_\alpha(y) &\leq C \int_{0}^{2^{k+1}R} \T_{|x|}^\alpha \tilde{f}_e(|y|)  d\mu_\alpha(y)  \nonumber
    \\
    & \leq C (2^{k+1}r)^{d_\alpha} (\tilde{f}_e)^*(|x|) \nonumber
    \end{align}
    where $(\tilde{f}_e)^*$ is the maximal function of $\tilde{f}_e$ on the Bessel–Kingman hypergroups \cite{Del}. Then substituting the above estimate in \eqref{eq: 4.12a}, we obtain
    \begin{align}
    |E_{B^1(0,r)}(x)| &\leq C  (\tilde{f}_e)^*(|x|) \sum_{k= -\infty}^{1} \rho(2^k r) \\
% \end{align}
% Since 
% \begin{align}
%     \int_{2^k r}^{2^{k+1}r} \frac{\rho(t)}{t} dt \geq C \rho(2^k r)  \int_{2^k r}^{2^{k+1}r} \frac{1}{t} dt = C \rho(2^k r) \ln 2, \nonumber
% \end{align}
% using \eqref{eq:3.12} and the above estimate in \eqref{eq:4.12}, we have
% \begin{align}
     &\leq C (\tilde{f}_e)^*(|x|) \sum_{k=-\infty}^{1} \int_{2^k r}^{2^{k+1}r} \frac{\rho(t)}{t} dt \nonumber
    \\
    &\leq C (\tilde{f}_e)^*(|x|) \int_0^{4r} \frac{\rho(t)}{t} dt \nonumber
    \\ 
    &\leq C (\tilde{f}_e)^*(|x|) \int_0^{r} \frac{\rho(t)}{t} dt\leq C \frac{\psi(r)}{\phi(r)} \tilde{f}_e^*(|x|) , \label{eq:4.13}
\end{align}
where we have used \eqref{eq:3.12}. \\
For $E_{B^2(0,r)}$, we have
\begin{align}
    |E_{B^2(0,r)}(x)| \leq \int_{B^c(0,2r)} |\tau_{x_0}^\alpha f(y) | \left| \tau_x^\alpha \frac{\rho(|y|)}{|y|^{d_\alpha}} - \frac{\rho(|y|)}{|y|^{d_\alpha}} \right| d\mu_\alpha(y). \nonumber
\end{align}
Then Lemma \ref{lem: l7.1} and Lemma \ref{lem:4.1} with $k=1$ imply that $E_{B^2(0,r)}$ is finite and 
\begin{align}
  |E_{B^2(0,r)}(x)| &\leq C r \int_{B^c(0,2r)} |\tau_{x_0}^\alpha f(y)| \frac{\rho(|y|)}{|y|^{d_\alpha +1}} d\mu_\alpha(y)  \nonumber\\
      & = C r \sum_{k = 1}^{\infty} \int_{2^k r \leq |y| < 2^{k+1}r} \frac{\rho(|y|)}{|y|^{d_\alpha +1}}  |\tau_{x_0}^\alpha f(y)|  d\mu_\alpha(y) \nonumber
    \\
    & \leq C r\sum_{k= 1}^{\infty} \frac{\rho(2^k r)}{( 2^k r)^{d_\alpha +1}}\int_{2^k r \leq |y| < 2^{k+1}r}  |\tau_{x_0}^\alpha f(y)|  d\mu_\alpha(y). \label{eq:47a}
    \end{align}
Since the operator $\tau_x^\alpha$ is linear it suffices to estimate the above integral for non-negative functions only. Therefore using Proposition \ref{prop:4.1b}, for each $k=1,2,3,...$ we can write 
\begin{align}
    \int_{2^k R \leq |y| < 2^{k+1} R} |\tau_{x_0}^\alpha f(y)| d\mu_\alpha(y) & \leq C \int_{0}^{2^{k+1}R} \T_{|x_0|}^\alpha f_e(|y|)  d\mu_\alpha(y)  \nonumber
    \\
    & = C\int_0^\infty f_e(y) \T_{|x_0|}^\alpha \chi_{J(0, 2^{k+1}R)}(y) d\mu_\alpha(y) \nonumber
    \\
      & = C\int_0^\infty f_e(y) (\T_{|x_0|}^\alpha \chi_{J(0, 2^{k+1}R)}(y))^{\frac{1}{p} + \frac{1}{p'}} d\mu_\alpha(y), \nonumber
\end{align}
where $p'$ is conjugate exponent of $p> 1$, i.e. $\frac{1}{p} + \frac{1}{p'} =1$. Applying H\"older's inequality,  we get
\begin{align}
    \int_{2^k R \leq |y| < 2^{k+1} R} |\tau_{x_0}^\alpha f(y)| d\mu_\alpha(y) 
    & \leq C \left( \int_0^\infty (f_e)^p(y) \T_{|x_0|}^\alpha \chi_{J(0, 2^{k+1} R)}(y) d\mu_\alpha(y) \right)^{\frac{1}{p}}\nonumber\\
    &\hspace{1.75cm}\times\left( \int_0^\infty \T_{|x_0|}^\alpha \chi_{J(0, 2^{k+1} R)}(y)d\mu_\alpha(y) \right)^{\frac{1}{p'}}. \nonumber 
\end{align}
Now  using the support property and the boundedness of the generalized translation $\T_{|x|}^\alpha$ on Bessel–Kingman
hypergroups (see Proposition \ref{prop:1}), we obtain
\begin{align}
  &\int_{2^k R \leq |y| < 2^{k+1} R } |\tau_{x_0}^\alpha f(y)| d\mu_\alpha(y)\nonumber\\
  &\leq C \left( \int_{J(|x_0|, 2^{k+1} R)} (f_e)^p(y)\T_{|x_0|}^\alpha \chi_{J(0, 2^{k+1} R)} d\mu_\alpha(y) \right)^{\frac{1}{p}} \left( \int_0^{2^{k+1} R}  d\mu_\alpha(y) \right)^{\frac{1}{p'}} \nonumber
  \\
   &\leq C \left( \int_{J(|x_0|, 2^{k+1} R)} (f_e)^p(y) d\mu_\alpha(y) \right)^{\frac{1}{p}} ( 2^{k+1}R)^{\frac{1}{p'}}. \nonumber
   \end{align}
 As we have done earlier, in the above step, using the same observation, we obtain
   \begin{align}
   \int_{2^k R \leq |y| < 2^{k+1} R } |\tau_{x_0}^\alpha f(y)| d\mu_\alpha(y) & \leq C  \left( \int_{I(|x_0|, 2^{k+1} R)} (f_e)^p(y) d\mu_\alpha(y) \right)^{\frac{1}{p}} ( 2^{k+1}R)^{\frac{1}{p'}} \nonumber
    \\
    & = C  \left( \int_{B(0, 2^{k+1} R)} \tau_{|x_0|}^\alpha (f_e)^p(y) d\mu_\alpha(y) \right)^{\frac{1}{p}} ( 2^{k+1}R)^{\frac{1}{p'}}\nonumber
    \\
     & \leq C \phi(2^k R) (2^k R)^{d_\alpha} ||f_e||_{L^{p,\phi}(\mathbb{R}, d\mu_\alpha)} \nonumber
     \\
    & \leq C \phi(2^k R) (2^k R)^{d_\alpha} ||f||_{L^{p,\phi}(\mathbb{R}, d\mu_\alpha)}. 
\end{align}
Proceeding similarly as in the proof of Lemma \ref{lem:4.1} with $j=1$, we obtain from \eqref{eq:47a}
\begin{align}
  |E_{B^2(0,r)}(x)|  & \leq C r||f||_{L^{p, \phi}(\mathbb{R}, d\mu_\alpha)}\sum_{k= 1}^{\infty} \frac{\rho(2^k r)}{ 2^k r} \phi(2^{k+1} r)  \nonumber
     \\
      &\leq C r ||f||_{L^{p, \phi}(\R, d\mu_\alpha)} \int_{r}^{ \infty} \frac{\rho(t) \phi(t)}{t^2} dt. \nonumber
\end{align}
Further using \eqref{eq:3.12} once again, we can write
\begin{align}
   |E_{B^2(0,r)}(x)| \leq C \psi(r) ||f||_{L^{p, \phi}(\R, d\mu_\alpha)}. \label{eq:4.14}
\end{align}
Now for $1 < p < \infty$, it follows from \eqref{eq:51c}
\begin{align}
    |\tau_{x_0}^\alpha \tilde{T}_\rho^\alpha f(x) - C_{B(0,r)}(x_0)|^p = |E_{B(0,r)}(x)|^p \leq 2^{p-1} (|E_{B^1(0,r)}(x)|^p + |E_{B^2(0,r)}(x)|^p). \nonumber
\end{align}
  Hence using  [\cite{bloom}, p. 58] (see also \cite{pak}) the boundedness of the maximal operator on the Bessel–Kingman hypergroups and using \eqref{eq:4.13}, we get
\begin{align}
    \frac{1}{\psi(r)^p r^{d_\alpha}} \int_{B(0,r)} |E_{B^1(0,r)}(x)|^p d\mu_\alpha(x) &\leq \frac{1}{\phi(r)^p r^{d_\alpha}} \int_{B(0,r)} [(\tilde{f}_e)^*(|x|)]^p d\mu_\alpha(x) \nonumber \\
    &\leq C  \frac{1}{\phi(r)^p r^{d_\alpha}}  || (\tilde{f}_e)^*||_{L^{p}(\R, d\mu_\alpha)}^p \nonumber
    \\
    &\leq C  \frac{1}{\phi(r)^p r^{d_\alpha}}  || \tilde{f}_e||_{L^{p}(\R, d\mu_\alpha)}^p \nonumber
    \\
    &\leq C  \frac{1}{\phi(r)^p r^{d_\alpha}} ||\tilde{f}||^p_{L^{p}(\R, d\mu_\alpha)} \leq C ||f||^p_{L^{p,\phi}(\R, d\mu_\alpha)} , \nonumber
\end{align}
and using \eqref{eq:4.14}, we get
\begin{align}
    \frac{1}{\psi(r)^p r^{d_\alpha}} \int_{B(0,r)} |E_{B^2(0,r)}(x)|^p d\mu_\alpha(x) \leq C ||f||^p_{L^{p,\phi}(\R, d\mu_\alpha)}. \nonumber
\end{align}
Finally combining these two estimates, we obtain
\begin{align}
    \frac{1}{\psi(r)^p r^{d_\alpha}} \int_{B(0,r)} |\tau_{x_0}^\alpha \tilde{T}_\rho^\alpha f(x) - C_{B(0,r)}(x_0)|^p d\mu_\alpha(x) \leq C_p ||f||^p_{L^{p,\phi}(\R, d\mu_\alpha)}, \nonumber
\end{align}
and hence the theorem follows from the definition of campanato norm \eqref{eq:3.8}.\qed

 \textbf{Proof of Theorem \ref{3.5}.} Let $f \in \mathcal{L}^{p,\phi}(\mathbb{R}, d\mu_\alpha)$. For given $r >0$, let $\tilde{B}=B(0,2r)$ and suppose $x \in B(0,r)$.
  Consider
  \begin{align}
     & E_{B(0,r)}(x) = \int_\mathbb{R} ( \tau_{x_0}^\alpha f(y) - f^\alpha_{\tilde{B}}(x_0) ) \left( \tau_x^\alpha \frac{\rho(|y|)}{|y|^{d_\alpha}} - \frac{\rho(|y|)}{|y|^{d_\alpha}} (1 - \chi_{\tilde{B}}(y))\right) d\mu_\alpha(y),  \nonumber
     \\
      & C_{B^1(0,r)}(x_0) = \int_\mathbb{R} ( \tau_{x_0}^\alpha f(y) - f^\alpha_{\tilde{B}}(x_0) ) \nonumber\\
      &\hspace{4cm} \times\left( \frac{\rho(|y|)}{|y|^{d_\alpha}}(1 - \chi_{\tilde{B}}(y)) - \frac{\rho(|y|)}{|y|^{d_\alpha}} (1 - \chi_{B_0}(y))\right) d\mu_\alpha(y), \nonumber
      \\
       & C_{B^2(0,r)}(x_0) = \int_\mathbb{R} f^\alpha_{\tilde{B}}(x_0)  \left( \tau_{x}^\alpha \frac{\rho(|y|)}{|y|^{d_\alpha}}  - \frac{\rho(|y|)}{|y|^{d_\alpha}} (1 - \chi_{B_0}(y))\right) d\mu_\alpha(y), \nonumber
       \\
         & E_{B^1(0,r)}(x) = \int_{\tilde{B}} ( \tau_{x_0}^\alpha f(y) - f^\alpha_{\tilde{B}}(x_0) ) \tau_x^\alpha \frac{\rho(|y|)}{|y|^{d_\alpha}} d\mu_\alpha(y), \hspace{6.25 cm} \nonumber
         \\
          & E_{B^2(0,r)}(x) = \int_{\tilde{B}^c} ( \tau_{x_0}^\alpha f(y) - f^\alpha_{\tilde{B}}(x_0) ) \left( \tau_x^\alpha \frac{\rho(|y|)}{|y|^{d_\alpha}} - \frac{\rho(|y|)}{|y|^{d_\alpha}}\right) d\mu_\alpha(y),  \nonumber 
  \end{align}
  where $x_0 \in \mathbb{R}$. Then from \eqref{eq:3.4}
  \begin{align}
      \tau_{x_0}^\alpha \tilde{T}_\rho^\alpha f(x) - (C_{B^1(0,r)}(x_0) + C_{B^2(0,r)}(x_0)) = E_{B(0,r)}(x) = E_{B^1(0,r)}(x) + E_{B^2(0,r)}(x). \nonumber
  \end{align}
As shown in \cite{spsa}, $C_{B^1(0,r)}(x_0)$ and $C_{B^2(0,r)}(x_0)$ are well-defined.
Now consider
\begin{align}
    E_{B^1(0,r)}(x) &= \int_{B(0,2r)} \left(\tau_{x_0}^\alpha f(y) - f^\alpha_{B(0,2r)}(x_0)\right) \tau_x^\alpha \frac{\rho(|y|)}{|y|^{d_\alpha}} d\mu_\alpha(y) \nonumber\\
    &= \int_{\R} \left((\tau_{x_0}^\alpha f - f^\alpha_{B(0,2r)}(x_0))\chi_{B(0,2r)}\right)(y) \tau_x^\alpha \frac{\rho(|y|)}{|y|^{d_\alpha}} d\mu_\alpha(y) \nonumber\\
    &= \int_{B(-x, 2r)} \tau_{-x}^\alpha \left((\tau_{x_0}^\alpha f - f^\alpha_{B(0,2r)}(x_0))\chi_{B(0,2r)}\right)(y) \frac{\rho(|y|)}{|y|^{d_\alpha}} d\mu_\alpha(y). \nonumber
\end{align}
We write $\tilde{f} := (\tau_{x_0}^\alpha f - f_{B(0,2r)}(x_0))\chi_{B(0,2r)}$ and $\tilde{\phi}(r) = \int_r^\infty \frac{\phi(t)}{t} dt$ to have
\begin{align}
    E_{B^1(0,r)}(x) = \int_{B(-x, 2r)} \tau_{-x}^\alpha \tilde{f}(y) \frac{\rho(|y|)}{|y|^{d_\alpha}} d\mu_\alpha(y). \nonumber
\end{align}
Now proceeding similarly as in the proof of Theorem \ref{3.4} and using \eqref{eq:3.13} , we obtain
\begin{align}
    |E_{B^1(0,r)}(x)| \leq C (\tilde{f}_e)^*(|x|) \int_0^{r} \frac{\rho(t)}{t} dt\leq C \frac{\psi(r)}{\tilde{\phi}(r)} (\tilde{f}_e)^*(|x|), \label{eq:4.13a}
\end{align}
where $(\tilde{f}_e)^*$ denotes the maximal function of $\tilde{f}_e$ on the Bessel–Kingman hypergroups \cite{Del}, and $\tilde{f}_e$ is the even part of $\tilde{f}$.
 Hence using [\cite{bloom}, p. 58] (see also \cite{pak}) the boundedness of the maximal operator on the Bessel–Kingman hypergroups, we get
\begin{align}
    &\frac{1}{\psi(r)}\left( \frac{1}{r^{d_\alpha}} \int_{B(0,r)} |E_{B^1(0,r)}(x)|^p d\mu_\alpha(x) \right)^{\frac{1}{p}} \nonumber
    \\
    &\leq \frac{C}{\tilde{\phi}(r) r^{d_\alpha/p}} \left( \int_{B(0,r)} [ (\tilde{f}_e)^*(|x|)]^p d\mu_\alpha(x) \right)^{1/p} \nonumber
    \\
    & \leq \frac{C}{\tilde{\phi}(r) r^{d_\alpha/p}} ||\tilde{f}||_{L^p(\R, d\mu_\alpha)} \nonumber
     \\
     % & \leq \frac{C}{\tilde{\phi}(r) r^{d_\alpha/p}}  ||\tilde{f} - \sigma (f)_{\chi_{B(0,2r)}} ||
    & \leq \frac{C}{\tilde{\phi}(r) r^{d_\alpha/p}}  ||\tilde{f} - \sigma (f)_{\chi_{B(0,2r)}} + \sigma (f)_{\chi_{B(0,2r)}}||_{L^p(\R,d\mu_\alpha)} \nonumber
    \\
    & \leq \frac{C}{\tilde{\phi}(r) r^{d_\alpha/p}} \Bigg( ||(\tau_{x_0}^\alpha f - \sigma(f) )_{\chi_{B(0,2r)}}||_{L^p(\R,d\mu_\alpha)} \nonumber
    \\
    &\quad +||(f_{B(0,2r)}(x_0) - \sigma (f) )_{\chi_{B(0,2r)}} ||_{L^p(\R, d\mu_\alpha)} \Bigg) \nonumber 
    \\
     % &= \frac{C}{\tilde{\phi}(r) r^{d_\alpha/p}} \Bigg( 
     % \|(\tau_{x_0}^\alpha f - \sigma(f) )_{\chi_{B(0,2r)}}\|_{L^p(\mathbb{R}, d\mu_\alpha)} \nonumber
     % \\
     % &\quad + \left( \int_{\mathbb{R}} |f_{B(0,2r)}(x_0) - \sigma(f)|^p \chi_{B(0,2r)}(y) \, d\mu_\alpha(y) \right)^{1/p} 
     % \Bigg) \nonumber
     % \\
    &  = \frac{C}{\tilde{\phi}(r) r^{d_\alpha/p}}  \Bigg( \left(\int_{B(0,2r)}|(\tau_{x_0}^\alpha f - \sigma(f) )|^p(y) d\mu_\alpha(y) \right)^{1/p} \nonumber
    \\
    &\quad +  |f^\alpha_{B(0,2r)}(x_0) - \sigma(f)| r^{d_\alpha/p} \Bigg) \nonumber 
    \\
      &  = \frac{C}{\tilde{\phi}(r) r^{d_\alpha/p}}  \left( \int_{B(0,2r)}|\tau_{x_0}^\alpha (f - \sigma(f) )|^p(y) d\mu_\alpha(y) \right)^{1/p} \nonumber
      \\
      &\quad +  |f^\alpha_{B(0,2r)}(x_0) - \sigma(f)| r^{d_\alpha/p} \label{eq: 4.23c}
    \end{align}
   where $\sigma(f) = \lim_{r \to \infty} f_{B(0,r)}(x_0)$. Now first we will estimate 
    \begin{align}
       \left( \int_{B(0,2r)}|\tau_{x_0}^\alpha (f - \sigma(f) )|^p(y) d\mu_\alpha(y) \right)^{1/p}. \nonumber 
    \end{align}
   Since the operator $\tau_x^\alpha$ is linear it suffices to estimate the above integral for non-negative functions only. Let  $f^\#
$ denote the non-negative part of the function $f - \sigma(f)$. Then, using Proposition \ref{prop:4.1b}, we get
    \begin{align}
       \left( \int_{B(0,2r)}|\tau_{x_0}^\alpha f^\#(y)|^p d\mu_\alpha(y) \right)^{1/p}  &\leq C_p \left( \int_{0}^{2r} \T_{|x_0|} (f^\#_e)^p(y) d\mu_\alpha(y) \right)^{1/p} \nonumber
       \\
       & = C_p \left( \int_{0}^{\infty} \T_{|x_0|} (f^\#_e)^p(y) \chi_{J(0, 2r)}(y) d\mu_\alpha(y) \right)^{1/p} \nonumber
       \\
        & = C_p \left( \int_{0}^{\infty}  (f^\#_e)^p(y) \T_{|x_0|}\chi_{J(0, 2r)}(y) d\mu_\alpha(y) \right)^{1/p}. \nonumber
    \end{align}
    Using the support property of the translation on the Bessel-Kingman hypergroups given in Proposition \ref{prop:1}, it follows at once that
    \begin{align}
         \left( \int_{B(0,2r)}|\tau_{x_0}^\alpha f^\#(y)|^p d\mu_\alpha(y) \right)^{1/p}  & \leq C_p \left( \int_{J(|x_0|, 2r)}  (f^\#_e)^p(y) \T_{|x_0|}\chi_{J(0, 2r)}(y) d\mu_\alpha(y) \right)^{1/p} \nonumber
         \\
         &\leq C_p \left( \int_{J(|x_0|,2r)}  (f^\#_e)^p(y) d\mu_\alpha(y) \right)^{1/p}. \nonumber
    \end{align}
  As we have done previously, we apply the same observation in the above step and consequently obtain the following
    \begin{align}
         \left( \int_{B(0,2r)}|\tau_{x_0}^\alpha f^\#(y)|^p d\mu_\alpha(y) \right)^{1/p}
         &\leq C \left( \int_{I(|x_0|,2r)}  (f^\#_e)^p(y)  d\mu(y) \right)^{1/p} \nonumber
         \\
         & = C \left( \int_{B(0,2r)} \tau_{|x_0|}^\alpha (f^\#_e)^p(y)  d\mu_\alpha(y) \right)^{1/p} \nonumber
         \\
         &\leq C \tilde{\phi}(2r)(2r)^{d_\alpha/p} ||f^\#_e||_{L^{p, \tilde{\phi}}(\R, d\mu_\alpha)} \nonumber
          \\
           &\leq C \tilde{\phi}(r)r^{d_\alpha/p} ||f^\#||_{L^{p, \tilde{\phi}}(\R, d\mu_\alpha)}. \label{eq: 4.26} 
    \end{align}
   By putting the value from \eqref{eq: 4.26} in \eqref{eq: 4.23c}, we obtain
    \begin{align}
    \frac{1}{\psi(r)}&\left( \frac{1}{r^{d_\alpha}} \int_{B(0,r)} |E_{B^1(0,r)}|^p d\mu_\alpha(x) \right)^{\frac{1}{p}} \nonumber
    \\
    &\leq \frac{C}{\tilde{\phi(r)} r^{d_\alpha/p}} \left[ ||f- \sigma(f) ||_{L^{p,\tilde{\phi}}(\R, d\mu_\alpha)} \tilde{\phi}(r) r^{d_\alpha/p} + |f^\alpha_{B(0,2r)} - \sigma(f)| r^{d_\alpha/p} \right] \nonumber
    \\
    & = C \left[ ||f- \sigma(f)||_{L^{p,\tilde{\phi}}(\R, d\mu_\alpha)} + \frac{1}{\tilde{\phi}(r)} |f^\alpha_{B(0,2r)}(x_0) - \sigma(f)| \right]. \label{eq:57c}
\end{align}
Moreover, Lemma \ref{lem:4.3} and Lemma \ref{lem:7} lead to the following inequalities 
\begin{align}
    ||f- \sigma(f)||_{L^{p,\tilde{\phi}}(\R, d\mu_\alpha)} \leq C ||f||_{\mathcal{L}^{p, \phi}(\R, d\mu_\alpha)}, \nonumber
\end{align}
and
\begin{align}
  |f_{B(0,2r)}(x_0) - \sigma(f)| \leq C_2 ||f||_{\mathcal{L}^{p,\phi}(\R, d\mu_\alpha)} \tilde{\phi}(r). \nonumber  
\end{align}
Finally, using these inequalities we otain from \eqref{eq:57c}
\begin{align}
    \frac{1}{\psi(r)}\left( \frac{1}{r^{d_\alpha}} \int_{B(0,r)} |E_{B^1(0,r)}(x)|^p d\mu_\alpha(x) \right)^{\frac{1}{p}} \leq C ||f||_{\mathcal{L}^{p,\phi}(\R, d\mu_\alpha)}.  \label{eq: 4.28}
\end{align}
For $E_{B^2(0,r)}$, already we have
\begin{align}
    |E_{B^2(0,r)}(x)| &\leq \int_{B^c(0,2r)} |\tau_{x_0}^\alpha f(y) - f^\alpha_{B(0,2r)}(x_0) | \left| \tau_x^\alpha \frac{\rho(|y|)}{|y|^{d_\alpha}} - \frac{\rho(|y|)}{|y|^{d_\alpha}} \right| d\mu_\alpha(y).\nonumber
\end{align}
Then Lemma \ref{lem: l7.1} implies that 
\begin{align}
  |E_{B^2(0,r)}(x)| &\leq C r \int_{B^c(0,2r)} |\tau_{x_0}^\alpha f(y) - f^\alpha_{B(0, 2r)}(x_0)| \frac{\rho(|y|)}{|y|^{d_\alpha +1}} d\mu_\alpha(y)  \nonumber\\
      & = C r \sum_{k = 1}^{\infty} \int_{2^k r \leq |y| < 2^{k+1}r} \frac{\rho(|y|)}{|y|^{d_\alpha +1}}  |\tau_{x_0}^\alpha f(y) - f^\alpha_{B(0, 2r)}(x_0)|  d\mu_\alpha(y) \nonumber
    \\
    & \leq C r\sum_{k= 1}^{\infty} \frac{\rho(2^{k+1} r)}{( 2^{k+1} r)^{d_\alpha +1}}\int_{2^k r \leq |y| < 2^{k+1}r}  |\tau_{x_0}^\alpha f(y) - f^\alpha_{B(0, 2r)}(x_0)|  d\mu_\alpha(y) \nonumber
    \\
     & \leq C r\sum_{k= 1}^{\infty} \frac{\rho(2^{k+1} r)}{( 2^{k+1} r)^{d_\alpha +1}} \nonumber \\
     &\hspace{0.75cm}\times\left(\int_{2^k r \leq |y| < 2^{k+1}r}  |\tau_{x_0}^\alpha f(y) - f^\alpha_{B(0, 2r)}(x_0)|^{p}  d\mu_\alpha(y) \right)^{1/p} (2^{k+1} r)^{d_\alpha(1-\frac{1}{p})}  \nonumber
     \\
     & \leq C r \sum_{k=1}^\infty \frac{\rho(2^{k+1} r)}{2^{k+1} r} \nonumber\\
     &\hspace{2.25cm}\times\left( \frac{1}{(2^k r)^{d_\alpha}} \int_{0 \leq |y| < 2^{k+1}r} |\tau^\alpha_{x_0} f(y) - f^\alpha_{B(0, 2r)}(x_0)|^{p} d\mu_\alpha(y) \right)^{1/p}. \label{eq: 4.29} 
    \end{align}
Using the Lemma \ref{lem:4.2} in the above \eqref{eq: 4.29}, we get
\begin{align}
  |E_{B^2(0,r)}(x)|  & \leq C r ||f||_{\mathcal{L}^{p, \phi}(\mathbb{R}, d\mu_\alpha)}\sum_{k= 1}^{\infty} \frac{\rho(2^{k+1} r)}{ 2^{k+1} r} \int_{2r}^{2^{k+1} r} \frac{\phi(t)}{t} dt  \nonumber
     \\
      & \leq C r ||f||_{\mathcal{L}^{p, \phi}(\mathbb{R}, d\mu_\alpha)}\sum_{k= 1}^{\infty} \int_{2^{k+1} r}^{2^{k+2}r } \frac{\rho(t)}{t^2} \left( \int_{2r}^t \frac{\phi(s)}{s} ds \right) dt  \nonumber
     \\
     &\leq C r ||f||_{\mathcal{L}^{p, \phi}(\R, d\mu_\alpha)} \int_{2 r}^{\infty} \frac{\rho(t)}{t^2} \left( \int_{2r}^t \frac{\phi(s)}{s} ds \right) dt \nonumber
     \\
     &\leq Cr ||f||_{\mathcal{L}^{p, \phi}(\R, d\mu_\alpha)} \int_{2r}^\infty \left( \int_s^{\infty} \frac{\rho(t)}{t^2} dt \right) \frac{\phi(s)}{s} ds. \nonumber
     \end{align}
     Using \eqref{eq:3.5} and then \eqref{eq:3.13}, it implies that
    \begin{align}
     |E_{B^2(0,r)}(x)| \leq C r ||f||_{L^{p, \phi}(\R, d\mu_\alpha)} \int_{2r}^{ \infty} \frac{\rho(s) \phi(s)}{s^2} ds  \leq C \psi(r) ||f||_{\mathcal{L}^{p, \phi}(\R, d\mu_\alpha)}. \nonumber
\end{align}
This follows that
\begin{align}
    \frac{1}{\psi(r)} \left( \frac{1}{\mu_\alpha(B(0,r))} \int_{B(0,r)} |E_{B(0,r)}(x)|^p d\mu_\alpha(x) \right)^{1/p} \leq C ||f||_{\mathcal{L}^{p,\phi}(\R, d\mu_\alpha)}. \label{eq:4.30}
\end{align}
Summing the estimates \eqref{eq: 4.28} and \eqref{eq:4.30}, we obtain
\begin{align}
    &\frac{1}{\psi(r)} \left( \frac{1}{\mu_\alpha(B(0,r))} \int_{B(0,r)} |\tau_{x_0}^\alpha \tilde{T}_\rho^\alpha f(x) - (C_{B^1(0,r)(x_0)} + C_{B^2(0,r)}(x_0))|^p d\mu_\alpha(x) \right)^{1/p} \nonumber
    \\
   & \hspace{11cm} \leq C_p ||f||_{\mathcal{L}^{p,\phi}(\R, d\mu_\alpha)}. \nonumber
\end{align}
This completes the proof.\qed

\begin{remark}
    Note that the authors proved the boundedness of $\tilde{T}_\rho^\alpha$ from  $\mathcal{L}^{1, \phi}(\R,d\mu_\alpha)$ to $\mathcal{L}^{1,\psi}(\R, d\mu_\alpha)$ in \cite{spsa} (see Theorem $7.1$) using the Fubini's theorem. However Fubini's theorem can't be applicable for the case $1<p<\infty$, therefore in this article, we provide the proof of the boundedness of $\tilde{T}_\rho^\alpha$ on the generalized Dunkl-type Campanato space when $1 < p < \infty$ using the boundedness  of  the maximal operator on the Bessel–Kingman hypergroups. 
\end{remark}

%%%%%%%%%%%%%%%%%%%%%%%%%% BACK MATTERS %%%%%%%%%%%%%%%%%%%%%%%%%%%%
 \section*{
  Declarations} %%%%%%%%%%%%%%%%%%%%
\textbf{Competing Interests} The author declares no competing interests\\

\textbf{Ethical Approval} 
    Not applicable  \\
    
 \textbf{Funding}
First author is supported by the University Grants Commission(UGC) with Fellowship No. 211610060698/(CSIRNETJUNE2021). Second author is supported by  seed grant IITM/SG/SWA/94 from Indian Institute of Technology Mandi, India. \\

\textbf {Data availability statement} Not applicable\vspace{3mm}
%%%%%%%%%%%%%%%%%% REFERENCES: %%%%%%%%%%%%%%%%%%%%%%%%%%%%%%%%%%%%%%%%%%%%
%% BibTeX users: please use \bibliographystyle{spmpsci} %% for math. and phys. sci.
%% Non-BibTeX users: please use the model as below !!! %%

%%%% for FCAA - pls. include directly the Refs items here ! %%%
%%%% following STRICTLY the models below %%%%%
%%%% and ARRANGE the items in ALPHABETIC ORDER for authors' family names !!!

 %%%%%%%%%%%%%%%%%%%%%%%%%%%%%%

\bigskip  %%%%%%%%%%%%%%%%%%%%%%%%%%%%%%%%%

\small %%%
\noindent
{\bf Publisher's Note}
Springer Nature remains neutral with regard to jurisdictional claims in published maps and institutional affiliations. 
\end{document}